\documentclass[a4paper,12pt]{amsart}

\usepackage{graphicx}
\usepackage{amsmath}
\usepackage{amssymb}
\usepackage{amscd}
\input xy
\xyoption{all}

\newcommand{\comment}[1]{}

\newcommand{\ko}{\: , \;}
\newcommand{\ul}[1]{\underline{#1}}

\renewcommand{\tilde}[1]{\widetilde{#1}}

\renewcommand{\top}{\operatorname{top}\nolimits}
\renewcommand{\mod}{\operatorname{mod}\nolimits}

\newcommand{\rad}{\operatorname{rad}\nolimits}
\newcommand{\soc}{\operatorname{soc}\nolimits}

\newcommand{\idim}{\operatorname{id}\nolimits}

\newcommand{\End}{\operatorname{End}\nolimits}

\newcommand{\gdim}{\operatorname{gl.dim}\nolimits}

\newcommand{\Ext}{\operatorname{Ext}\nolimits}
\newcommand{\Hom}{\operatorname{Hom}\nolimits}
\newcommand{\C}{\mathcal{C}}
\newcommand{\N}{\mathbb{N}}
\newcommand{\Z}{\mathbb{Z}}

\newcommand{\repL}{{\widehat \Lambda}}

\newcommand{\modL}{\textrm{mod} \Lambda}

\newcommand{\modRepL}{\textrm{mod} \repL}
\newcommand{\stmodRepL}{\underline{\textrm{mod}} \repL}
\newcommand{\Ainf}{\mathbf{A}_\infty}
\newcommand{\ZAinf}{\Z\Ainf}
\newcommand{\citeARS}{[ARS]}
\newcommand{\citeF}{[F]}
\newcommand{\citeHI}{[H1]}
\newcommand{\citeHII}{[H2]}
\newcommand{\citeR}{[R]}

\DeclareMathOperator{\Img}{Im}

\newcommand{\secI}{Section~1}
\newcommand{\secII}{Section~2}
\newcommand{\secIII}{Section~3}

\renewcommand{\Im}{\operatorname{Im}\nolimits}
\newcommand{\Ker}{\operatorname{Ker}\nolimits}
\newcommand{\Coker}{\operatorname{Coker}\nolimits}

\newcommand{\pd}{\operatorname{pd}\nolimits}
\newcommand{\uline}{\underline}
\newcommand{\la}{\Lambda}
\newcommand{\mcal}{\mathcal}
\newcommand{\xto}{\xrightarrow}
\newcommand{\what}{\widehat}
\newcommand{\wtilde}{\widetilde}

\newcommand{\vphi}{\varphi}

\newtheorem{theorem}{Theorem}[section]

\newtheorem{corollary}[theorem]{Corollary}

\newtheorem{example}[theorem]{Example}

\newtheorem{lemma}[theorem]{Lemma}

\newtheorem{proposition}[theorem]{Proposition}

\newcommand{\arr}[1]{\stackrel{#1}{\longrightarrow}}
\newcommand{\iso}{\stackrel{_\sim}{\rightarrow}}
\newcommand{\id}{\mathbf{1}}
\newcommand{\can}{\mbox{{\rm can}}}
\newcommand{\cc}{{\mathcal C}}

\newcommand{\ce}{{\mathcal E}}

\newcommand{\ci}{{\mathcal I}}

\newcommand{\cm}{{\mathcal M}}

\newcommand{\cp}{{\mathcal P}}
\newcommand{\cR}{{\mathcal R}}

\newcommand{\ten}{\otimes}
\newcommand{\lten}{\overset{\boldmath{L}}{\ten}}

\newcommand{\p}{\mathbf{p}}

\begin{document}
\title{Bounded derived categories and repetitive algebras}
\author{Dieter Happel, Bernhard Keller and Idun Reiten}
\maketitle
\section*{Introduction}
Let $\la$ be a finite dimensional algebra over a field $k$. It was
proved in
\cite{H1} that there is a full and faithful embedding of the bounded
derived
category  ${\bf D}^b(\la)$ into the stable category
$\uline{\mod}\what{\la}$
of finite dimensional modules over the repetitive algebra
$\what{\la}$. This
embedding is an equivalence if and only if $\la$ has finite global
dimension
\cite{H1}. The category ${\bf D}^b(\la)$ is a triangulated category
which
does not have almost split triangles when $\la$ has infinite  global
dimension, whereas $\uline{\mod}\what{\la}$ is triangulated and always has
almost
split triangles.

The purpose of this paper is to investigate the relationship between
${\bf D}^b(\la)$ and $\uline{\mod}\what{\la}$ from various points of
view,
which is of course meaningful only for algebras $\la$ of infinite
global
dimension. The most satisfactory results are obtained for Gorenstein
algebras,
especially for selfinjective algebras.

We investigate the embedding ${\bf D}^b(\la)\subset
\uline{\mod}\what{\la}$
from  the point of view of universal properties with respect to
triangle
functors to triangulated categories with almost split triangles, and
also to
which extent $\uline{\mod}\what{\la}$ is the  smallest category
containing
${\bf D}^b(\la)$ with these properties. The first question has a
positive
answer for Gorenstein algebras, and is not true in general. The second
question has a negative answer even for selfinjective algebras.

We also investigate the behavior of almost split triangles and
irreducible
maps under the embedding functor, and show that both are actually
preserved.
While it is known from \cite{H2} what the end terms of almost split
triangles
in ${\bf D}^b(\la)$ look like, and hence the left and right end terms
of certain
irreducible maps, we do not know in general so much about irreducible
maps
in ${\bf D}^b(\la)$. However in the selfinjective case we show that
there are
no irreducible maps  not associated with almost split triangles when
$(\rad\la)^2\ne 0$, and we describe them all when $(\rad\la)^2=0$. We
believe
that also for arbitrary $\la$ there should be very few irreducible
maps  not
associated with almost split triangles. As an application  of our
results we
show that when $\la$ is selfinjective all the components of the
AR-quiver of
the category ${\bf K}^b(\mcal{P})$ of bounded  complexes of projective
modules
are of the form $\mathbb{Z}A_{\infty}$.

The paper is organised as follows.
In section 1 we give some background material from \cite{H1} on the
categories
${\bf D}^b(\la)$ and $\mod\what{\la}$, including properties of almost
split
triangles. In section 2, we give an example
showing that in general, the embedding ${\bf D}^b(\la) \subset
\uline{\mod}\what{\la}$ is not universal among triangle functors
from ${\bf D}^b(\la)$
to triangulated categories with almost split triangles.
We also show that the embedding
$\mod\la\subset \uline{\mod}\what{\la}$ has a weak universal
property with
respect to triangle functors from ${\bf D}^b(\la)$ to triangulated
categories where the Nakayama functor becomes an equivalence.
We deduce that if $\la$ is Gorenstein, there is a natural
triangle functor from $\uline{\mod}\what{\la}$ to ${\bf D}^b(\la)$.
In section 3 we show that
even
when  $\la$ is selfinjective, there is  an infinite strictly descending chain
of
triangulated subcategories of $\uline{\mod}\what{\la}$ with almost
split
triangles and containing ${\bf D}^b(\la)$. In section 4 we show that
irreducible maps in ${\bf D}^b(\la)$  stay irreducible in
$\uline{\mod}\what{\la}$ , and give sufficient conditions for the
existence of
irreducible maps in ${\bf D}^b(\la)$ of  the form  $S[-1]\to T$  where
$S$ and
$T$ are simple $\la$-modules. In section 5  we show that almost split
triangles in ${\bf D}^b(\la)$ stay almost split in
$\uline{\mod}\what{\la}$,
and give  the shape of the components of the AR-quiver of
${\bf K}^b(\mcal{P})$ for
selfinjective
algebras. We also give  necessary conditions for having  irreducible
maps not
coming from almost split triangles for Gorenstein algebras, and deduce
the
result on irreducible maps in ${\bf D}^b(\la)$ when $\la$ is
selfinjective.
In section 6 we deal with arbitrary finite dimensional algebras, and
give
some results supporting the suspicion that there are very few
irreducible
maps not associated with almost split triangles. We also show that the
natural
questions of a connection between irreducible maps between infinite
complexes
of projective modules and their finite parts have negative answers.
\section{Preliminaries}
 In this section we will fix the notation and recall some of the results
 frequently
used in the subsequent sections. For the proofs of the stated
propositions we
refer to \cite{H1}. Let $\la$ be a finite dimensional algebra over a
field $k$.

 We denote by $\mod\la$ the category of finitely generated left
 $\la$-modules
and by $_{\la}\mcal{P}$ (resp. $_{\la}\mcal{I}$) the full subcategory
of projective (resp. injective) $\la$-modules. For a simple
$\la$-module $S$
we denote by $P(S)$ (resp. $I(S)$ ) the projective cover
(resp. injective
envelope) of $S$. We denote by
$\nu_{\la}\colon_{\la}{\mcal{P}}\to _{\la}\mcal{I}$ the Nakayama
functor defined
by $\nu_{\la}=D\Hom(-,_{\la}\la)$, where $D$ is the duality with
respect to $k$,
and by $\nu^-_{\la}\colon _{\la}\mcal{I}\to _{\la}\mcal{P}$ the
inverse Nakayama
functor which is defined by $\nu^-_{\la}=\Hom(D\la_{\la},-)$. We
denote by
${\bf D}^b(\la)$ the bounded derived category of $\mod\la$. The
Nakayama
functors $\nu_{\la}$ and $\nu^-_{\la}$ induce inverse equivalences of
triangulated
 categories still denoted by
$\nu_{\la}\colon {\bf K}^b(_{\la}\mcal{P})\to {\bf
  K}^b(_{\la}\mcal{I})$ and
$\nu^-_{\la}\colon {\bf K}^b(_{\la}{\mcal{I}})\to {\bf
  K}^b(_{\la}\mcal{P})$.
We denote by ${\bf K}^{-,b}(_{\la}\mcal{P})$ the homotopy category of
complexes
over $_{\la}\mcal{P}$ bounded above with bounded cohomology groups.
Note that ${\bf K}^{-,b}(_{\la}\mcal{P})\simeq {\bf D}^b(\la)$. For a
complex $Z=(Z^i,d^i)$ in ${\bf D}^b(\la)$ and $n\in \mathbb{Z}$ we
always have
a triangle $Z_{\geq n}\to Z\to Z_{<n}\to Z_{\geq n}[1]$ in ${\bf
  D}^b(\la)$,
where $Z_{\geq n}^i=Z^i$ for $i\geq n$, $Z_{<n}=Z^i$ for $i<n$, and
zero otherwise,
with the induced differentials.

To $\la$ we may associate the repetitive algebra $\what{\la}$ and its
category
$\mod\what{\la}$ of finitely generated modules. The
$\what{\la}$-modules $X$
are given by $X=(X_i,f_i)$ where $X_i\in\mod\la$ and $X_i=0$ for
almost all $i$,
 $f_i\colon X_i\to \nu^-_{\la}X_{i+1}$ such that
 $f_i\nu^-_{\la}(f_{i+1})=0$ for all $i$.
A morphism of $\what{\la}$-modules is defined in an obvious way. There
is an
automorphism $\nu_{\what{\la}}\colon \mod\what{\la}\to \mod\what{\la}$
defined by
  $(\nu_{\what{\la}}X)_i=X_{i+1}$. The inverse is denoted by
  $\nu^-_{\what{\la}}$.
The category $\mod\what{\la}$ is a Frobenius category in the sense of
\cite{H1}. The indecomposable projective-injective
$\what{\la}$-modules are
given by $I=P=(X_i,f_i)$ with $X_i=P(S), X_{i+1}=I(S),f_i=\idim_{P(S)}$
and zero
otherwise. Note that $\top P=\nu_{\what{\la}}\soc P$. Clearly there
are enough
projectives. So for each $X$  we obtain exact sequences
$0\to X\to I(X)\to \Omega^-_{\what{\la}}X\to 0$ and
$0\to \Omega_{\what{\la}}X\to P(X)\to X\to 0$. We denote by
$\uline{\mod}\what{\la}$ the stable category. This is a triangulated
category
where $\Omega^-_{\what{\la}}$ serves as a translation functor. If
$X\in \uline{\mod}\what{\la}$, we may choose a representative again
denoted
by $X\in\mod\what{\la}$ without indecomposable projective direct
summands.
This fact will be used frequently later on.

There is a triangle functor
$\mu\colon{\bf D}^b(\la)\to \uline{\mod}\what{\la}$ which is full and
faithful such that $\mu$ extends the identity functor on $\mod\la$
where $\mod\la$ is embedded in ${\bf D}^b(\la)$ (resp.
$\uline{\mod}\what{\la}$) as complexes (resp. modules) concentrated in
degree zero. It is known \cite{H2} that $\mu$ is an equivalence if and
only
if $\gdim \la < \infty$.

In general we recall from \cite{GK} the following
description
 of $\Im \mu$. $Z=(Z_i,g_i)\in \Im\mu$ if and only if  there is some
 $n\ge 0$
such that $(\Omega^{-n}_{\what{\la}}Z)_j=0$  for $j>0$ and
$(\Omega_{\what{\la}}^nZ)_j=0$
 for $j<0$. Also note that for a $\what{\la}$-module $Z=(Z_i,g_i)$
 with $Z_i=0$
for $i<0$ also $(\Omega^r_{\what{\la}}Z)_j=0$ for $j<0$ and all $r\ge 0$.

This has the following immediate consequence for Gorenstein algebras (see also \cite{CZ}).

\begin{corollary} Let $\Lambda$ be a Gorenstein algebra. Then
  $\text{im}\,\mu=\{Z=(Z_i,g_i)\in\text{mod}\,\what{\la}\,|\, \text{pd}_\Lambda
  Z_i<\infty\,\,\text{for}\,\,i\neq 0\}$
\end{corollary}
\begin{proof} If $\Lambda$ is a Gorenstein algebra then the modules of finite
  projective dimension coincide with the modules of finite injective dimension.
Moreover this dimension is bounded by the projective dimension of
$D\Lambda_\Lambda$ which coincides with the injective dimension of
${}_\Lambda\Lambda.$ Suppose that
$Z=(Z_i,g_i)\in\text{mod}\,\what{\la}$ satisfies $\text{pd}_\Lambda
  Z_i<\infty$ for $i\neq 0,$ then it follows immediately from the criterion
  mentioned above from \cite{GK} that $Z\in\text{im}\,\mu.$ Conversely let
$Z=(Z_i,g_i)\in\text{im}\,\mu$ and assume that
$Z=(Z_i,g_i)\in\text{mod}\,\what{\la}$ satisfies ${pd}_\Lambda Z_i=\infty$
for some $i\neq 0.$ We may assume that $i>0.$ Choose $i$ maximal with this
property. So $\text{pd}_\Lambda Z_j<\infty$ for $j>i.$ By the first part of
the proof and the fact that $\text{im}\,\mu$
is a triangulated category the factor module $Z'=(Z'_j,g_j)$ with $Z_j'=Z_j$
for $j\le i$ and $Z_j'=0$ for $j>i$ is contained in $\text{im}\,\mu.$
But then $(\Omega^{-n}_{\what{\la}}Z')_i\neq 0$ for all $n\ge 0,$
in contrast to \cite{GK}
\end{proof}
Let $\mcal{C}$ be  a triangulated category which is Krull-Schmidt. Let
$Z$
be an object in $\mcal{C}$. We say that there is an almost split
triangle
ending at $Z$  provided there is a triangle in $\mcal{C}$ of the form
$$X\xto{u}Y\xto{v}Z\xto{w}X[1]$$
 where (i) $X$ is indecomposable, (ii) for all $f\colon W\to Z$ not split
 epi
there is some $g\colon W\to Y$ with $f=gv$ and (iii) $w\ne 0$.

We refer to \cite{H1} for equivalent formulations and the connection
to
irreducible maps.

In case there is  an almost split triangle ending at $Z$, the starting
term
$X$ is uniquely determined up to isomorphism. We then define $\tau_{\mcal{C}}Z=X$.

It easily follows from the existence of almost split sequences in
$\mod\what{\la}$ that $\uline{\mod}\what{\la}$ has almost split
triangles.
It is well known and can be shown using the definition of $\tau$ that
for
$Z\in \uline{\mod}\what{\la}$ indecomposable, then
$\tau_{\what{\la}}Z=\nu_{\what{\la}}\Omega^2_{\what{\la}}X$.

In the case of ${\bf D}^b(\la)$ the following is known \cite{H2}. Let
$Z\in{\bf D}^b(\la)$ be indecomposable. Then there is an almost split
triangle $X\to Y\to Z\xto{w}X[1]$ if and only if
$Z\in {\bf K}^b(_{\la}\mcal{P})$. In this case $\tau_{{\bf  D}^b(\la)}Z=\nu_{\la} Z[-1]$.
 Thus  ${\bf D}^b(\la)$ has almost split triangles if and only if
$\gdim \la< \infty$.

\section{A counterexample and a weak universal property}

\noindent {\bf Problem:} Let $\Lambda$ be an artin algebra and
\[
\mu : \mathbf{D}^b(\Lambda) \to \ul{\mod}\, \widehat{\Lambda}
\]
the embedding of \cite{H1}. Let $\cc$ be a triangulated category
with Auslander-Reiten triangles and $F: \mathbf{D}^b(\Lambda) \to
\cc$ a triangle functor. Does there exist a triangle functor
\[
G: \ul{\mod}\, \widehat{\Lambda} \to \cc
\]
such that $F \iso \mu G$ ?

\bigskip\noindent
The following example shows that the answer is no, in general.

\bigskip\noindent
{\bf Example:} Let $\Lambda$ be given as a factor algebra of a path
algebra of a field $k$ by an ideal:
\[
\xymatrix{T \ar@(ul,dl)_\beta & S \ar[l]^\alpha } \ko \langle
\beta^2, \alpha\beta\rangle.
\]
Let $S,T$ be the two simple $\Lambda$-modules. Then
\[
P(S) = \left( \begin{array}{c} S \\ T \end{array} \right) \mbox{ and
} P(T) = \left( \begin{array}{c} T \\ T \end{array} \right)
\]
are the indecomposable projective $\Lambda$-modules and $I(S)=S$ and
\[
I(T)=\left( \begin{array}{cc} S \hspace*{0.5cm} T \\ T \end{array}
\right)
\]
the indecomposable injective $\Lambda$-modules. It is easy to see
that $\Lambda$ is not Gorenstein ($I(T)$ is of infinite projective
dimension). Let $\Gamma=\mbox{End}_\Lambda(P(T))$, so
$\Gamma=k[x]/(x^2)$ and let $F=\Hom(P(T), - )$ be a functor from
$\mod \Lambda$ to $\mod \Gamma$. Now $F$ is exact, so $F$ induces a
functor $\mathbf{D}^b(\Lambda)\to \mathbf{D}^b(\Gamma)$. Since
$\Gamma$ is selfinjective, there is a functor \cite{ric}
\[
\pi: \mathbf{D}^b(\Gamma) \to \ul{\mod}\, \Gamma \ko
\]
so there is a triangle functor $\phi: \mathbf{D}^b(\Lambda) \to
\ul{\mod}\, \Gamma$ and $\ul{\mod}\, \Gamma$ has Auslander-Reiten
triangles. We are now going to show that there is no triangle
functor $G: \ul{\mod}\,\widehat{\Lambda} \to \ul{\mod}\, \Gamma$
such that $\phi=\mu G$.

Suppose there exists a triangle functor $G:
\ul{\mod}\,\widehat{\Lambda} \to \ul{\mod}\, \Gamma$ such that
$\phi=\mu G$. Let $X=(X_i, f_i)$ be an object of
$\ul{\mod}\,\widehat{\Lambda}$ with $X_1=S$ and $X_i=0$ for $i\neq
1$. Then $\Omega^-_{\widehat{\Lambda}} X =P(S)$, the stalk module
concentrated in degree zero. So $G\Omega^-_{\widehat{\Lambda}} X
=\phi(P(S))=T$ and
\[
G\Omega^-_{\widehat{\Lambda}} X \cong \Omega^-_\Gamma G(X) = GX \ko
\]
so $GX\cong T$. Also $G(S)=\phi(S)=0$ and $G(T)=\phi(T)=T$. But then
also $G(\Omega_{\widehat{\Lambda}} T)=T$. Now
$\Omega_{\widehat{\Lambda}} T=(Y_i, f_i)$ where $Y_0=T$, $Y_1=I(T)$
and $f_0: T \to P(T)$ the canonical map, and $Y_i=0$ for $i\neq
0,1$. Consider the exact sequence in $\mod \Lambda$:
\[
0 \to T \to P(T) \to T \to 0.
\]
It gives rise to a triangle
\[
T[-1] \arr{f} T \to P(T) \to T \quad (*)
\]
in $\mathbf{D}^b(\Lambda)$. Since $\phi(P(T))=0$, the map $\phi(f)$
is invertible. Now we also have the exact sequence of $\mod
\widehat{\Lambda}$:
\[
0 \to Z \to \Omega_{\widehat{\Lambda}} T \to T \to 0.
\]
This gives rise to a triangle
\[
Z \to \Omega_{\widehat{\Lambda}} T \arr{\mu(f)} T \to Z[1] \ko
\]
which identifies with the image of the triangle $(*)$ under $\mu$.
Applying $G$ then shows that $GZ=0$ because $\mu G(f)=\phi(f)$ is
invertible. Now $Z=(Z_i, f_i)$, where $Z_1=I(T)$ and $Z_i=0$ for
$i\neq 1$. Let $U=(U_i, f_i)$ with $U_1=T$ and $U_i=0$ for $i\neq
1$. Then we obtain an exact sequence in $\mod \widehat{\Lambda}$
\[
0 \to U \to Z \to U\oplus X \to 0 \ko
\]
which gives rise to a triangle
\[
U \to Z \to U\oplus X \to U[1]
\]
in $\ul{\mod}\, \widehat{\Lambda}$ and so
\[
GU \to GZ \to GU\oplus GX \to GU[1]
\]
is a triangle in $\ul{\mod}\, \Gamma$. Now $GZ=0$ by the computation
above and $G(X)=T$, so the triangle is of the form
\[
GU \to 0 \to GU\oplus T \to GU[1] \ko
\]
a contradiction.

\bigskip
\noindent {\bf A weak universal property.} As the above
counterexample shows, the repetitive category is not the `universal
triangulated category with Auslander-Reiten triangles containing the
derived category'. However, we will see that if we take into account
additional structure, we do get a weak universal property for the
embedding
\[
\mod\Lambda \to \ul{\mod}\,\hat{\Lambda}.
\]
Roughly speaking this embedding is the `universal functor to a
triangulated category where the Nakayama functor becomes an
equivalence'. In the case where $\Lambda$ is Gorenstein, we will use
this property to construct a natural triangle functor from the
stable category of the repetitive category to the bounded derived
category.

Let us now construct the additional structure we need: For short,
let us write $\cm$ for $\mod\Lambda$ and $\cR$ for
$\mod\hat{\Lambda}$. We write $\Sigma :\cm\to\cm$ for the right
exact extension of the Nakayama functor defined in section~1: Thus,
we have $\Sigma M = (D\Lambda)\ten_\Lambda M$ for all $M$ in $\cm$.
We now define an exact functor $\cR\to\cR$, which we will also
denote by $\Sigma$. Namely, we put
\[
\Sigma X = \nu_{\hat{\Lambda}} (\Omega X) \ko
\]
where $\Omega$ is the syzygy functor $\cR\to\cR$ constructed as
follows: If $X$ is an object of $\cR$ with structure maps $f_i$,
$i\in\Z$, we define the object $PX$ to have the $i$th component
\[
(\Lambda \ten_k X_i) \oplus (D\Lambda \ten_k X_{i-1})
\]
and the structure maps
\[
\left[ \begin{array}{cc} 0 & 0 \\ \id & 0 \end{array} \right] :
\nu(PX)_i \to (PX)_{i-1}.
\]
Thus, the object $PX$ is projective-injective. We define the
canonical map $PX \to X$ to have the components
\[
[\can, g_{i-1}] : (PX)_i \to X_i
\]
where $\can$ is the canonical map from $\Lambda \ten_k X_i$ to $X_i$
and $g_{i-1}$ is the map $D\Lambda\ten_k X_{i-1} \to X_i$ induced by
$\nu (f_{i-1})$. Thus, the map $PX \to X$ is a functorial projective
right approximation of $X$. We define $\Omega X$ to be the kernel of
$PX \to X$.

The functor $\Sigma : \cR \to \cR$ is exact, preserves
projective-injectives and induces an equivalence in the stable
category (namely, the Serre functor). Moreover, if $F_0 : \cm\to\cR$
denotes the canonical embedding, we have a morphism of functors
\[
\phi_0 : F_0 \Sigma \to \Sigma F_0.
\]
Namely, for an object $M$ of $\cm$, the only non vanishing component
of the morphism $F_0\Sigma (M) \to  \Sigma F_0(M)$ is induced by the
canonical map $D\Lambda \ten_k M \to D\Lambda \ten_\Lambda M$. It is
easy to check that if $P$ is a projective $\Lambda$-module, then
$\phi_0 P$ becomes an isomorphism in the stable category of $\cR$.
To summarize, we have
\begin{itemize}
\item a $k$-linear Frobenius category $\cR$ endowed with an exact functor $\Sigma:\cR\to\cR$
preserving projective-injectives and inducing an equivalence in the
stable category,
\item an exact functor $F_0 : \cm \to \cR$ endowed with a morphism
\[
\phi_0 : F_0\Sigma \to  \Sigma F_0
\]
such that $\phi_0 P$ becomes an isomorphism in the stable category
for each projective module $P$.
\end{itemize}

\begin{theorem} Let $\ce$ be a $k$-linear Frobenius category endowed with
an exact functor $\Sigma : \ce \to \ce$ preserving
projective-injectives and inducing an equivalence in the stable
category. Let $F: \cm \to \ce$ be an exact functor endowed with a
morphism $\phi: \Sigma F \to F \Sigma$ such that $\phi P$ becomes an
isomorphism in the stable category for each projective module $P$.
Then there is a triangle functor
\[
G : \ul{\cR} \to \ul{\ce} \ko
\]
such that $G$ commutes with $\Sigma$ up to isomorphism and the
triangle
\[
\xymatrix{\cm \ar[r]^{F_0} \ar[dr]_F & \ul{\cR} \ar[d]^G \\
 & \ul{\ce}
}
\]
commutes up to isomorphism. \comment{ More precisely, there is an
isomorphism of triangle functors $\psi: G\Sigma  \to  \Sigma G$ and
an isomorphism of functors $\gamma : F_0 G\to F$ such that we have
\[
 \phi \circ \Sigma\gamma = (\gamma \Sigma)(G\phi_0)(\psi F_0).
\]
\[
\xymatrix{ \Sigma G F_0 \ar[d]_{\Sigma \gamma} \ar[r]^{\psi F_0} &
G\Sigma F_0 \ar[r]^{G\phi_0} &
GF_0 \Sigma \ar[d]^{\gamma \Sigma} \\
\Sigma F \ar[rr]_{\phi} & & F\Sigma }
\]
}
\end{theorem}

The theorem will be proved below. Note that it does not make any
claim about uniqueness. In fact, one could obtain a more intrinsic
formulation and a uniqueness statement by working in a more
sophisticated framework based on towers of triangulated categories
\cite{k2}, or derivators \cite{g} or the homotopy category of dg
categories \cite{ta} \cite{to} \cite{k1}. However, this would go
beyond the scope of this article.

\begin{corollary} Suppose that $\Lambda$ is Gorenstein. Then there
is a a triangle functor
\[
G: \ul{\mod}\,\widehat{\Lambda} \to \mathbf{D}^b(\Lambda)
\]
which commutes with the inclusion of $\mod\Lambda$ and such that we
have a functorial isomorphism
\[
\Sigma \circ G \iso G \circ \Sigma  \ko
\]
where $\Sigma : \ul{\mod}\,\widehat{\Lambda} \to
\ul{\mod}\,\widehat{\Lambda}$ is the Serre functor and $\Sigma :
\mathbf{D}^b(\Lambda) \to \mathbf{D}^b(\Lambda)$ the functor $M
\mapsto D\Lambda \lten_\Lambda M$.
\end{corollary}

\bigskip
\noindent {\em Proof of the corollary.} Let $\ce$ be the category of
right bounded complexes of projective $\Lambda$-modules with bounded
homology. Then the stable category of $\ce$ is triangle equivalent
to the bounded derived category. For each $\Lambda$-bimodule $B$,
write $\p(B)$ for a projective bimodule resolution of $B$. Let
$\Sigma : \ce \to \ce$ be the (total) tensor product over $\Lambda$
by the complex of bimodules $\p(D\Lambda)$. Let $F_0$ be the functor
taking a module $M$ to $\p(\Lambda)\ten_\Lambda M$. To construct
$\phi: F_0 \Sigma\to \Sigma F_0$, it suffices to construct a
quasi-isomorphism of bimodule complexes
\[
\tilde{\phi} : \p(D\Lambda) \ten_\Lambda \p(\Lambda)  \to
\p(\Lambda) \ten_\Lambda D\Lambda.
\]
Indeed, since the morphism
\[
\p(D\Lambda) \ten_\Lambda \p(\Lambda) \to D\Lambda
\]
is a projective resolution, it lifts (in the homotopy category)
along the quasi-isomorphism
\[
\p(\Lambda) \ten_\Lambda D\Lambda \to D\Lambda
\]
and we define $\tilde{\phi}$ to be a representative of a lift. If
$P$ is a projective module, then in the square (commutative in the
homotopy category),
\[
\xymatrix{
\p(D\Lambda) \ten_\Lambda \p(\Lambda)\ten_\Lambda P \ar[r] \ar[d]_{\phi P} & D\Lambda\ten_\Lambda P \ar[d]^{\id} \\
\p(\Lambda) \ten_\Lambda D\Lambda \ten_\Lambda P \ar[r] &
D\Lambda\ten_\Lambda P }
\]
the two horizontal morphisms are quasi-isomorphisms and so the left
vertical morphism is a homotopy equivalence. This means that $\phi
P$ becomes an isomorphism in the stable category. Thus, the
hypotheses of the theorem are satisfied and we get, if $\Lambda$ is
Gorenstein, a natural triangle functor
\[
G : \ul{\mod}\,\widehat{\Lambda} \to \mathbf{D}^b(\Lambda)
\]
which extends the inclusion of $\mod\Lambda$ and commutes with
$\Sigma$ up to isomorphism of triangle functors.

\bigskip
\noindent {\em Proof of the theorem.} It is not hard to see that it
suffices to define a functor with the required properties on the
full subcategory of objects $X$ of $\cR$ with $X_i=0$ for $i>0$. Let
$X$ be such an object of $\cR$ with structure maps $f_i: \Sigma X_i
\to X_{i+1}$, $i\in\Z$. We define $G_1 X$ to be the complex over
$\ce$ with components $\Sigma^i F(X_{-i})$ and with the differential
\[
\Sigma^i (FX_{-i}) \to \Sigma^{i-1} (F X_{-i+1})
\]
given by $(\Sigma^{i-1}\phi X_{-i})(\Sigma^{i-1} Ff_i)$. It is
straightforward to check that the square of the differential
vanishes and that with the natural definition of $G_1$ on morphisms,
we get a $k$-linear functor
\[
G_1:\cR \to \mathbf{C}^b(\ce)
\]
taking exact sequences of $\cR$ to componentwise conflations of the
category $\mathbf{C}^b(\ce)$ of bounded complexes over $\ce$.
Moreover, the functor $G_1$ takes an indecomposable projective
injective object given by a projective $P$ (put in degree $-1$, for
simplicity of notation) and the identity $\Sigma P \to \Sigma P$ to
a complex of the shape
\[
\xymatrix{ \ldots \ar[r] & 0 \ar[r] & \Sigma (FP) \ar[r]^{\phi P}  &
F(\Sigma P) \ar[r] & 0 \ar[r] & \ldots }
\]
Now since $\ce$ is a Frobenius category, we have a canonical
triangle functor \cite{kv} \cite{ric}
\[
\mathbf{D}^b(\ce) \to \ul{\ce}
\]
extending the natural projection functor $\ce \to \ul{\ce}$. We
define $G_2$ to be the composition
\[
\mathbf{C}^b(\ce) \to \mathbf{D}^b(\ce) \to \ul{\ce}
\]
and we put $G_3=G_1\circ G_2 : \cR \to \ul{\ce}$. Then $G_3$ takes
projective-injectives to zero-objects: Indeed, a complex of the form
\[
\xymatrix{ \ldots \ar[r] & 0 \ar[r] & \Sigma (FP) \ar[r]^{\phi P}  &
F(\Sigma P) \ar[r] & 0 \ar[r] & \ldots }
\]
is the cone over the morphism $\phi P$ (between complexes
concentrated in degree $0$). Since $\phi P$ becomes invertible in
$\ul{\ce}$ by assumption, the image of the cone under $G_2$ is a
zero object. Thus, $G_3$ induces a $k$-linear functor $G$. It is
clear from the construction that $F_0G$ is isomorphic to $F$. Since
$G_1$ takes conflations to componentwise conflations and the
projection $\mathbf{C}^b(\ce) \to \mathbf{D}^b(\ce)$ transforms each
componentwise conflation into a canonical triangle, the functor $G$
is in fact a triangle functor. Therefore, to construct a commutation
isomorphism $\Sigma G\to G\Sigma$ it suffices to construct such a
commutation isomorphism for $\Omega^{-1}\circ \Sigma$. Now in $\cR$,
the composition $\Omega^{-1}\circ \Sigma$ is isomorphic to the
degree shifting functor $\nu_{\widehat{\Lambda}}$. For an object
$X$, the image $G_1 (\nu_{\widehat{\Lambda}}X)$ is isomorphic to
$\Sigma (G_1 X) [-1]$, where we denote by $\Sigma$ the functor from
$\mathbf{C}^b(\ce)$ to itself obtained by applying
$\Sigma:\ce\to\ce$ to each component. Now the canonical triangle
functor
\[
\mathbf{D}^b(\ce) \to \ul{\ce}
\]
is functorial with respect to exact functors preserving
projective-injectives and thus canonically commutes with $\Sigma$.
Moreover, since it is a triangle functor, it is compatible with
shifts. So we get a canonical isomorphism
\[
\Omega^{-1}\Sigma G (X) \iso G_2 ( \Sigma (G_1 X)[-1]) \iso G_1 G_2
\Sigma \Omega^{-1}(X).
\]
\

Note that for Gorenstein algebras we now have a positive answer
to the question posed in the beginning of the section.

\section{Infinite chain of subcategories}
In this section we construct triangulated subcategories of
$\uline{\mod}\widehat{\la}$ containing $\mcal{D}^b(\la)$ for $\la$ a
selfinjective algebra. Recall from \secI{} that in this case
$\mcal{D}^b(\la)$ can be identified with the full subcategory of
$\uline{\mod}\widehat{\la}$ with objects $X=(X_i,f_i)$ such that $X_i$
is a
projective $\la$-module for $i\ne 0$.
\begin{lemma}\label{3.1}
  Let $\la$ be a selfinjective algebra. Let $I\subseteq\mathbb{Z}$ and
  let
$\mcal{C}_{I}\subset \{X=(X_i,f_i)\in\uline{\mod}\widehat{\la}\mid
X_i$ is a
projective $\la$-module if $i\not\in I\}$. Then the following hold

i) $\uline{\mcal{C}}_{I}\subseteq \uline{\mod}\widehat{\la}$ is a
triangulated
subcategory.

ii) $0\in I $ if and only if $\mcal{D}^b(\la) \subseteq \uline{\mcal{C}}_{I}$.

iii) If $I, I'\subseteq \mathbb{Z}$ with $I\subseteq I'$, then
$\uline{\mcal{C}_{I}}\subseteq \uline{\mcal{C}}_{I'}$.
\end{lemma}
\begin{proof}
  If $P$ is a projective-injective $\widehat{\la}$ -module, then
$P\in \mcal{C}_{I}$, since $\la$ is selfinjective. So
$\uline{\mcal{C}}_{I}\subseteq \uline{\mod}\widehat{\la}$.

Let $X\in \uline{\mcal{C}}_{I}$ and consider an exact sequence
$0\to X\to I(X)\to \Omega^-_{\widehat{\la}}X\to 0$ with  $I(X)$
injective in
$\mod\widehat{\la}$. Then for each $i\in\mathbb{Z}$ we have an exact
sequence
$0\to X_i\to I(X_i)\oplus \nu^- I(X_{i+1})\to Z_i\to 0$ where $I(X_i)$ is
the
$\la$-injective envelope of  $X_i$. If $i\not\in I$, the sequence
splits,
since $X_i$ is projective, hence $\Omega^-_{\widehat{\la}}X \in
\mcal{C}_{I}$.
Thus $\uline{\mcal{C}}_{I}$ is closed under the translation functor in
$\uline{\mod}\widehat{\la}$. Finally, let $X\to Y\to Z\to C[1]$ be a
triangle
in $\uline{\mod}\widehat{\la}$ with $X,Y\in
\uline{\mcal{C}_{I}}$. Then the
triangle gives an exact sequence $0\to X\to I(X)\oplus Y\to Z\to 0$ in
$\uline{\mod}\widehat{\la}$. So for each $j\in\mathbb{Z}$ we obtain an
exact
sequence $0\to X_i\to I(X)_i\oplus Y_i\to Z_i\to 0$  in $\mod\la$. If
$i\not\in I$ then $X_i,Y_i$ are projective, so $Z_i$ is projective,
hence
$Z\in\uline{\mcal{C}}_{I}$, so $\uline{\mcal{C}}_{I}$ is a
triangulated
subcategory of $\uline{\mod}\widehat{\la}$.

(ii) and (iii) are obvious.
\end{proof}
\begin{example}\label{3.2}
  Let $n\in\mathbb{N}$ and let $I=(n+1)\mathbb{Z}$. Let
$\mcal{C}_n=\mcal{C}_{I}$. Then $\nu^{n+1}_{\widehat{\la}}$ is an
automorphism on
$\mcal{C}_n$. In fact, if $X=(X_i,f_i)\in \mcal{C}_n$, then
$(\nu^{n+1}_{\widehat{\la}}X )_i=X_{i+n+1}$. Since $j\not\in
(n+1)\mathbb{Z}$ if and
only if $j+n+1\not\in(n+1)\mathbb{Z}$ we see that
$\nu^{n+1}_{\widehat{\la}}X\in \mcal{C}_n$.

If we choose $n+1=2^k$ and let $\mcal{D}_k=\mcal{C}_n$ we obtain a
descending
chain of subcategories
$\cdots \subseteq\uline{\mcal{D}}_2\subseteq\uline{\mcal{D}}_1
\subseteq\uline{\mod}\widehat{\la}$ and clearly
$\mcal{D}^b(\la)=\bigcap_{i\geq1}\uline{\mcal{D}}_i$.
\end{example}

Let $\la$ be a symmetric algebra and let
$F=\nu_{\widehat{\la}}\Omega_{\widehat{\la}}$
be the Serre functor on $\uline{\mod}\widehat{\la}$. So for all
$X,Y\in \mod\widehat{\la}$ we have
$\eta_{X,Y}\colon\uline{\Hom}(X,Y)\xrightarrow{\sim}D\uline{\Hom}(Y,FX)$
natural
in $X$ and $Y$. We will show that $F^{n+1}$ is a Serre functor on
$\uline{\mcal{C}}_n$. For this we will construct $\eta_{X}\colon
F^{n+1}X\to FX$
such that $\eta_{X}$ is natural in $X$ and for all $X,
Y\in\uline{\mcal{C}}_n$
we have that
$\uline{\Hom}(Y,F^{n+1}X)\xto{\sim}\uline{\Hom}(Y,FX)$. This then
implies that $F^{n+1}$ is a Serre functor on $\uline{\mcal{C}}_n$,
hence
$\uline{\mcal{C}}_n$ has Auslander-Reiten triangles.(Compare \cite{RV})
\begin{lemma}\label{3.3}
  For all $X\in\mod\widehat{\la}$ there is an exact sequence
$0\to K_{X}\xto{\mu_{X}}F(X)\xto{\pi_{X}}X\to 0$ which is natural in $X$.
\end{lemma}
\begin{proof}
  Since $\la$ is finite-dimensional, there is a functor
$P\colon \mod\la\to _{\la}\mcal{P}$ and an exact sequence
$0\to \Omega_{\la}X\xto{\alpha_{X}}P(X)\xto{\beta_{X}}X \to 0$ natural
in $X$
for $X\in \mod\la$ (compare section 2). Now let
$X=(X_i,f_i)\in\mod\widehat{\la}$. Applying
$\nu^-_{\widehat{\la}}$ if necessary we may assume that $X_i=0$ for
$i<0$. Now $P$
extends to a functor
$\widetilde{P}\colon\mod\widehat{\la}\to _{\widehat{\la}}\mcal{P}$ and
we have an
exact sequence $0\to\Omega_{\what{\la}}X\to\widetilde{P}(X)\to
X\to0$. Explicitly
we have for $i\geq0$ a commutative diagram of the form
$$\xymatrix@C2.5cm@R1cm{
    P(X_{i-1})\oplus \Omega_{\la}X_i\ar[r]^{
\left(
\begin{smallmatrix}
  P(f_{i-1})& 0\\
  \alpha_{X_i}& \Omega_{\la}(f_i)
\end{smallmatrix}\right)}\ar[d]_{
\left(
\begin{smallmatrix}
  1& 0\\
  0& \alpha_{X_i}
\end{smallmatrix}
\right)}& P(X_i)\oplus \Omega_{\la} X_{i+1}\ar[d]^{
\left(
\begin{smallmatrix}
  1& 0\\
  0& \alpha_{X_{i+1}}
\end{smallmatrix}\right)}\\
P(X_{i-1})\oplus P(X_i)\ar[r]^{
\left(
\begin{smallmatrix}
  -P(f_{i-1})& 0\\
  1& P(f_i)
\end{smallmatrix}\right)}\ar[d]_{
\left(
\begin{smallmatrix}
  o\\
  \beta_{X_i}
\end{smallmatrix}\right)}& P(X_i)\oplus P(X_{i+1})\ar[d]^{
\left(
\begin{smallmatrix}
  0\\
  \beta_{X_{i+1}}
\end{smallmatrix}\right)}\\
X_i\ar[r]^{f_i}& X_{i+1}
}$$
 The map $\pi_{X}\colon FX\to X$ is
 now defined by $(\pi_{X})_i\colon P(X_i)\oplus \Omega X_{i+1}\xto{\left(
\begin{smallmatrix}
\beta_{X_i}\\
0
\end{smallmatrix}\right)}X_i$. Clearly $\pi_{X}$  is surjective and
$K_{X}$ is
described by the following commutative  diagram
$$\xymatrix@C2.5cm@R1cm{
    \Omega(X_{i})\oplus \Omega_{\la}(X_{i+1})\ar[r]^{
\left(
\begin{smallmatrix}
  \Omega(f_{i})& 0\\
  1& \Omega_{\la}(f_{i+1})
\end{smallmatrix}\right)}\ar[d]_{
\left(
\begin{smallmatrix}
  \alpha_{X_i}& 0\\
  0& 1
\end{smallmatrix}
\right)}& \Omega(X_{i+1})\oplus \Omega (X_{i+2})\ar[d]^{
\left(
\begin{smallmatrix}
  \alpha_{X_{i+1}}& 0\\
  0& 1
\end{smallmatrix}\right)}\\
P(X_{i})\oplus \Omega(X_{i+1})\ar[r]^{
\left(
\begin{smallmatrix}
  P(f_{i})& 0\\
  \alpha_{X_i}& \Omega_{\la}(f_{i+1})
\end{smallmatrix}\right)}\ar[d]_{
\left(
\begin{smallmatrix}
  \beta_{X_i}\\
  0
\end{smallmatrix}\right)}& P(X_{i+1})\oplus \Omega(X_{i+2})\ar[d]^{
\left(
\begin{smallmatrix}
  \beta_{X_{i+1}}\\
  0
\end{smallmatrix}\right)}\\
X_i\ar[r]^{f_i}& X_{i+1}
}$$
Since  $P$, $\widetilde{P}$ are functors, the exact sequence
$0\to K_{X}\xto{\mu_{X}}F(X)\xto{\pi_{X}}X\to0$ is natural in $X$.
\end{proof}
So for each $1\le i\le n$ we obtain an exact sequence
$$(*)_i \text{   }0\to F^i K_{X}\to F^{i+1}X\to F^iX\to 0$$
Now $(*)_i$ induces $\eta_{X}\colon F^{n+1}X\to FX$ natural in $X$.
Thus for all
$X,Y\in\mod\what{\la}$ we have $\eta_{X,Y}\colon\Hom(Y,F^{n+1}X)\to
\Hom(Y,FX)$
natural in $X$ and $Y$.

For the following lemma we need some notation. Let $X\in \mod\la$ and
denote
by $\delta^i(X)=(Z_j,\gamma_j)$ the $\what{\la}$-module with
$Z_{i-1}=Z_i=X$, $Z_j=0$ for $j\ne i-1,i$, $\gamma_{i-1}=1_{X}$ and
$\gamma_j=0$
for $j\ne i-1$.
\begin{lemma}\label{3.4}
  If $X\in\mcal{C}_n$, then $K_{X}\simeq
  \bigoplus_{i\in(n+1)\mathbb{Z}}\delta^i
(\Omega X_i)$ in $\uline{\mod}\what{\la}$
\end{lemma}
\begin{proof}
  If $X=(X_i,f_i)\in\mcal{C}_n$, then by definition $X_i$ is
  projective for
$i\not\in(n+1)\mathbb{Z}$. Moreover it follows that $\delta^i(\Omega
X_i)$
is projective as a $\what{\la}$-module if $X_i$ is projective. It
follows
from the previous lemma that for each $i$ we have that
$\delta^i(\Omega X_i)$
is a submodule of $K_{X}.$ Explicitly consider the following commutative diagram
$$\xymatrix@C2cm{
    \Omega X_{i+1}\ar[r]^1\ar[d]^{(
      \begin{smallmatrix}
    0&1
      \end{smallmatrix})}& \Omega X_{i+1}\ar[d]^{(
      \begin{smallmatrix}
    1&\Omega f_{i+1}
      \end{smallmatrix})}\\
    \Omega X_i\oplus \Omega X_{i+1}\ar[r]^{\left(
      \begin{smallmatrix}
    \Omega f_i & 0\\
    1& \Omega f_{i+1}
    \end{smallmatrix}\right)}& \Omega X_{i+1}\oplus \Omega X_{i+2}
}$$
So if  $i\not\in(n+1)\mathbb{Z}$ we see that $\delta^i(\Omega X_i)$ is
a
direct summand of $K_{X}$. But then it follows from the description of
$K_{X}$
in the previous lemma that $K_{X}\simeq
\bigoplus_{i\in(n+1)\mathbb{Z}}\delta^i
(\Omega X_i)$ in $\uline{\mod}\what{\la}$
\end{proof}
\begin{lemma}\label{3.5}
  Let  $X\in \mod\la$ and $i\in\mathbb{Z}$. Then

i)$F\delta^i(X)\simeq\delta^{i-1}(\Omega X)$ in $\uline{\mod}\what{\la}$.

ii)If $Y\in\mcal{C}_n$, then $\uline{\Hom}(Y,\delta^i(X))=0$, if
$i\not\in(n+1)\mathbb{Z}$
\end{lemma}
\begin{proof}
  i) Clearly, if $X,Y\in\mod\la$ then
  $\Omega_{\what{\la}}\delta^i(X)=\delta^i
(\Omega_{\la} X)$ in $\uline{\mod}\what{\la}$ and
$\nu_{\what{\la}}\delta^i(Y)=
\delta^{i-1}(Y)$, so the assertion follows.

ii)Let $Y=(Y_i,g_i)\in\mcal{C}_n$ and let
$\vphi\in\Hom(Y,\delta^i(X))$. So we
have the following commutative diagram with $\vphi=(\vphi_j)$
$$\xymatrix{
  \cdots & Y_{i-2}\ar[r]^{g_{i-2}}\ar[d]&
  Y_{i-1}\ar[r]^{g_{i-1}}\ar[d]_{\vphi_{i-1}}&
Y_i\ar[r]^{g_i}\ar[d]^{\vphi_i}& Y_{i+1}\\
  \cdots & 0\ar[r]& X\ar[r]^{1_{X}}& X\ar[r]& 0
}$$
Consider $\delta^i(P(X))\xto{\pi}\delta^i(X)$, then $\delta^i(P(X))$
is a
projective $\what{\la}$-module. Since $i\not\in(n+1)\mathbb{Z}$, we
have
that $Y_i$ is a projective $\la$-module, so there is some
$\alpha_i\colon Y_i\to P(X)$ such that $\alpha_i\pi_i=\vphi_i$, where
$\pi=(\pi_j)$ and $\pi_j=0$ for $j\ne i-1,i$ and
$\pi_{i-1}=\pi=\beta_{X}$.

Set $\alpha_{i-1}=g_{i-1}\alpha_i$ and $\alpha_j=0$ for$j\ne
i-1,i$. Then
$\alpha=(\alpha_j)$ is a map $Y\to \delta^i(P(X))$ such that
$\alpha\pi=\vphi$,
 hence $\uline{\Hom}(Y,\delta^i(X))=0$
\end{proof}
\begin{proposition}\label{3.6}
  For all $X,Y\in \mcal{C}_n$, the natural transformation\\
  $\eta_{X,Y}\colon
\uline{\Hom}(Y,F^{n+1}X)\to \uline{\Hom}(Y,FX)$ is an isomorphism
\end{proposition}
\begin{proof}
  By the previous considerations we have for each $X\in\mcal{C}_n$ and
  each
$1\le i\le n$ an exact sequence $0\to F^iK_{X}\to F^{i+1}X\to
F^iX\to0$ in
$\mod\what{\la}$ which gives rise to a triangle $F^iK_{X}\to
F^{i+1}X\to F^iX\to
F^iK_{X}[1]$ in $\uline{\mod}\what{\la}$. Applying $\uline{\Hom}(Y,-)$
to this
triangle for $Y\in\mcal{C}_n$ gives  an exact sequence
$$\uline{\Hom}(Y,F^iK_{X}) \to \uline{\Hom}(Y,F^{i+1}X)\to
\uline{\Hom}(Y,F^iX)\to
 \uline{\Hom}(Y,F^iK_{X}[1])$$
By the description of $K_{X}$ and the previous lemma we see that
$\uline{\Hom}(Y,F^iK_{X})=0=\uline{\Hom}(Y,F^iK_{X}[1])$, hence
$\uline{\Hom}(Y,F^{i+1}K_{X}) \xto{\sim}\uline{\Hom}(Y,F^iX)$.

For each $1\le i\le n$ we also have an exact sequence
$0\to K^i\to F^{i+1}X\to FX\to 0$ in $\mod\what{\la}$ which gives rise
to a
triangle $K^i\to F^{i+1}X\to FX\to K^i[1]$  in
$\uline{\mod}\what{\la}$ where
$F^{i+1}X\to FX$ is obtained from the composition $F^{i+1}X\to F^iX\to
FX$. By the
 octahedral axiom we have a triangle $F^iK_{X}\to K^i\to K^{i-1}\to
 F^iK_{X}[1]$ .
By induction and the previous considerations we have that
$\uline{\Hom}(Y,K^i)=0=\uline{\Hom}(Y,K^i[1])$, hence
$\uline{\Hom}(Y,F^{i+1}X)\simeq\uline{\Hom}(Y,FX)$, thus $\eta_{X,Y}$
is an
isomorphism for all $X,Y\in\uline{\mcal{C}}_n$.
\end{proof}
As pointed out above this implies
\begin{corollary}\label{3.7}
  For each $n\in\mathbb{N}$ the category $\uline{\mcal{C}}_n$ has
  almost
split triangles
\end{corollary}

\section{Irreducible maps in ${\bf D}^b(\la)$}
In this section we study the behavior of the embedding
$\mu\colon{\bf D}^b(\la)\to \uline{\mod}\what{\la}$ under irreducible
maps.
If $X=(X^i,f^i)\in {\bf D}^b(\la)$ we write $\mu(X)=(X_i,f_i)$. If $X=(X^i,f^i)$
satisfies
 $X^i=0$ for $i<0$, then $X_i=0$ for $i<0$. Here  of course, we assume
 that
$\mu(X) $ has no projective-injective indecomposable summands. If
$X^i=0$ for
$i>0$ then $X_i=0$ for $i>0$.

We denote by $\mod^{\ge 0}\what{\la}=\{(X_i,f_i)|X_i=0, i<0\}$ and
$\mod^{>0}\what{\la}=\{(X_i,f_i)|X_i=0, i\le 0\}$. The categories
$\mod^{\le 0}\what{\la}$ and $\mod^{< 0}\what{\la}$ are defined
analogously.
Clearly $\mod^{\ge 0}\what{\la}$ and $\mod^{> 0}\what{\la}$ are stable
under
$\Omega_{\what{\la}}$. Moreover we clearly have $\Hom(X,Y)=0$ for
$X\in\mod^{\ge 0}\what{\la}$ and $Y\in \mod^{< 0}\what{\la}$. This
yields the
following easy lemma.
\begin{lemma}\label{4.1}
  Let $X\in \mod^{\ge 0}\what{\la}$ and $Y\in \mod^{<
    0}\what{\la}$. Then
$\Ext^1_{\what{\la}}(X,Y)=0$.
\end{lemma}
\begin{proof}
  We have $\Ext^1_{\what{\la}}(X,Y)\simeq
  \uline{\Hom}(X,\Omega^-_{\what{\la}}Y)=0$,
since $\Omega^-_{\what{\la}}Y\in\mod^{< 0}\what{\la}$.
\end{proof}
\begin{lemma}\label{4.2}
  Let $Z=(Z_i,f_i)\in\mod^{\ge 0}\what{\la}$ with $Z\in
  \Im\mu$. Consider the
exact sequence $0\to Z^{>0}\to Z\to Z_0\to 0$. Then
$\nu_{\what{\la}}Z^{>0}\in \Im\mu$.

\end{lemma}
\begin{proof}
  We verify the condition  mentioned in section 1 from \cite{GK}. Since
$\nu_{\what{\la}}Z^{>0}\in\mod^{\ge 0}\what{\la}$ we have that
$\Omega^r_{\what{\la}}(\nu_{\what{\la}}Z^{>0})\in\mod^{\ge
  0}\what{\la}$ for all
$r\ge 0$. The exact sequence $0\to Z^{>0}\to Z\to Z_0\to 0$ gives a
triangle
$\Omega_{\what{\la}}Z_0\to Z^{>0}\to Z\to Z_0$ in
$\uline{\mod}\what{\la}$. So for
each $n\ge 0$ we obtain a triangle
$$\Omega^{-n}_{\what{\la}}Z_0\to \Omega^{-n-1}_{\what{\la}}Z^{>0}\to
\Omega^{-n-1}_{\what{\la}}Z\to \Omega^{-n-1}_{\what{\la}}Z_0$$
For each $n\ge 0$ we clearly  have  $\Omega^{-n}_{\what{\la}}Z_0\in
\mod^{\le 0}\what{\la}$. Since $Z\in \Im\mu$ there is $n_0$ such that
$\Omega^{-n}_{\what{\la}}Z\in \mod^{\le 0}\what{\la}$ for all $n\ge
n_0$ by \cite{GK}.
Let $n\ge n_0$ and assume that
$\Omega^{-n-1}_{\what{\la}}Z^{>0}\not\in
\mod^{\le 0}\what{\la}$. Then there is $X\in \mod^{> 0}\what{\la}$
such that
$\uline{\Hom}(X,\Omega^{-n-1}_{\what{\la}}Z^{>0})\ne 0$. Since
$\uline{\Hom}(X,\Omega^{-n}_{\what{\la}}Z_0)=0=\uline{\Hom}(X,\Omega^{-n-1}_{\what{\la}}Z)$
 we obtain a contradiction.
\end{proof}
\begin{theorem}\label{4.3}
  Let $X,Y$ be indecomposable in ${\bf D}^b(\la)$ and let $f\colon
  X\to Y$ be
irreducible. Then $\mu(f)\colon\mu X\to \mu Y$ is irreducible in
$\uline{\mod}\what{\la}$.
\end{theorem}
\begin{proof}
We consider the almost split triangle
$\text{(*)   }\tau_{\what{\la}}\mu Y\xto{\alpha}E\xto{\beta}\mu Y\xto
{\gamma}\tau_{\what{\la}}\mu Y[1]$ in $\uline{\mod}\what{\la}$. Now
$\mu(f)$ is
not split epi, since $f$ is irreducible, hence we get $g\colon \mu
X\to E$
such that  $\mu(f)=g\beta$. Let $\mu X=(X_i,f_i)$ and $\mu
Y=(Y_i,g_i)$. We
may assume that $\mu X$, $\mu Y\in \mod^{\ge0}\what{\la}$. Now
$\tau_{\what{\la}}\mu Y=\nu_{\what{\la}}\Omega^2_{\what{\la}}\mu Y=(Z_i,h_i)$
satisfies $Z_i=0$
for $i<-1$, $Z_{-1}=\Omega^2_{\la}Y_0$. Thus $E=(E_i,u_i)$ satisfies
$E_i=0$ for
$i<-1$, $E_{-1}=\Omega^2_{\la}Y_0$ and $u_{-1}\colon \Omega^2_{\la}Y_0\to
\nu^-E_0$. Let
$P\xto{\pi}\Omega^2_{\la}Y_0$ be epi with $P$ a projective
$\la$-module. Let
$\widetilde{E}=(\widetilde{E}_i,v_i)\in \mod\what{\la}$ defined by
$\widetilde{E}_i=0$ for $i<-1$, $\wtilde{E}_{-1}=P$,
$\wtilde{E}_i=E_i$ for
$i\ge0$, $v_{-1}=\pi u_{-1}$, $v_i=u_i$ for $i\ge0$.

Now $\Omega^2_{\what{\la}}\mu Y\in \Im\mu$, since $\Im\mu$ is a
triangulated
subcategory. By Lemma \ref{3.2} we have that
$\nu_{\what{\la}}E^{>0}\in \Im\mu$.
Since $\nu_{\what{\la}}P\in\Im\mu$ we see that
$\wtilde{E}\in\Im\mu$. The
construction of $\wtilde{E}$ clearly yields a triangle
$$\nu_{\what{\la}}\mu(\Ker\pi)\xto{\delta}\wtilde{E}\xto{\epsilon}E\xto
{\eta}\nu_{\what{\la}}\mu(\Ker\pi)[1]$$
The factorization $\mu(f)=g\beta$ induces another factorization as follows:
$$\xymatrix{
  & \mu X\ar@{=}[r]\ar@{-->}[d]^{\wtilde{g}}& \mu X\ar[d]^g\\
  \nu_{\what{\la}}\mu(\Ker\pi)\ar[r]^{\delta}&
  \wtilde{E}\ar[r]^{\epsilon}\ar@{-->}
[d]^{\wtilde{\beta}}& E\ar[r]^{\eta}\ar[d]^{\beta}& \nu_{\what{\la}}\mu(\Ker\pi)[1]\\
  & \mu Y\ar@{=}[r]& \mu Y&
}$$
where $\wtilde{\beta}=\epsilon\beta$. Since $g\eta=0$ by Lemma
\ref{3.1} we
obtain $\wtilde{g}$ with  $\wtilde{g}\epsilon= g$. Now
$\wtilde{g}\wtilde{\beta}
=\wtilde{g}\epsilon\beta=g\beta=\mu(f)$.

Since $\wtilde{E}\in \Im \mu$ and $f$ is irreducible, we get that
$\wtilde{g}$
is a split mono or $\wtilde{\beta}$ is a split epi. If
$\wtilde{\beta}$ is
split epi, there is $\wtilde{\beta}_1\colon \mu Y\to \wtilde{E}$ such that
$\wtilde{\beta}_1\beta=1_{\mu Y}$. Set $\beta_1=\wtilde{\beta}_1c$. Then
$\beta_1\beta=\wtilde{\beta}_1c\beta=\wtilde{\beta}_1\wtilde{\beta}=1_{\mu Y}$,
so $\beta$ is a split epi, in contrast to (*)  being an almost split
triangle.
So $\wtilde{g}$ is split mono, hence there is some $\wtilde{g}_1\colon
\wtilde{E}\to \mu X$ such that $\wtilde{g}\wtilde{g}_1=1_{\mu X}$. Since
$\Hom(\nu_{\la}\Ker\pi,\mu X)=0$ , we have $\delta\wtilde{g}_1=0$, so
there is
some $g_1\colon E\to \mu X$ such that $\epsilon g_1=\wtilde{g}_1$. Now
$gg_1=\wtilde{g}\epsilon g_1=\wtilde{g}\wtilde{g}_1=1_{\mu X}$ shows that
$g$
is split mono, hence $X$ is an indecomposable direct summand of $E$
and
$\mu(f)$ is a component of $\beta$, hence $\mu(f)$ is irreducible.
\end{proof}
Next we show how certain irreducible maps in ${\bf D}^b(\la)$ arise
quite
naturally from extensions of simple $\la$-modules. This will be of
interest
in section 5.
\begin{proposition}\label{4.4}
  Let $S$ and $T$ be simple $\la$-modules with
  $\Ext_{\la}^1(S,T)\ne0$. If
$\rad P(S)$ and $ I(T)/T$ are both semisimple, then there is an
irreducible
map $f\colon S[-1]\to T$ in ${\bf D}^b(\la)$.
\end{proposition}
\begin{proof}
  We will show that there is an irreducible map $\vphi\colon
  \mu(S[-1])\to
\mu(T)$ in $\uline{\mod}\what{\la}$. For simplicity set $\mu(S)=S$ and
$\mu(T)=T$. Then $S[-1]=\Omega_{\what{\la}}S\simeq
\rad_{\what{\la}}P(S)$, where
$P(S)$ is the $\what{\la}$-projective cover of $S$. We consider the
almost
split sequence in $\mod\what{\la}$ starting in
$S[-1]=\rad_{\what{\la}}P(S)$. It
is well-known \cite{AR} that this is of the form
$$0\to\rad_{\what{\la}}P(S)\to
_{\what{\la}}P(S)\oplus\rad_{\what{\la}}P(S)/
\soc_{\what{\la}}P(S)\to P(S)/\soc_{\what{\la}}P(S)\to 0$$
Clearly, $\soc_{\what{\la}}P(S)=\nu^-_{\what{\la}}S$. Let $0\to
S\xto{\alpha}I(S)
\xto{\beta}I(S)/S\to0$ be exact in $\mod\la$, with $I(S)$ the
$\la$-injective
envelope of $S$. Applying the Nakayama functor
$\nu^{-}_{\la}=\Hom(D\la_{\la},-)$
yields an exact sequence
$$0\to
\Hom(D\la_{\la},S)\xto{\nu^-_{\la}(\alpha)}P(S)=\Hom(D\la_{\la},I(S))
\xto{\nu^-_{\la}( \beta)} \Hom(D\la_{\la},I(S)/S)$$
Let $g$ be the composition of $\rad P(S)\xto{\gamma}P(S)\xto
{\nu^-_{\la}( \beta)}\nu^-_{\la}I(S)/S$, so
$g=\gamma\nu^-_{\la}\beta$. Then
$\rad_{\what{\la}}P(S)/\soc_{\what{\la}}P(S)=(Z_i,g_i)=Z$ with
$Z_0=\rad_{\la}P(S)$,
$Z_1=I(S)/S$, $g_0=g$, and zero otherwise. By assumption, $\rad P(S)$
is
semisimple and $\Ext^1_{\la}(S,T)\ne0$, so $T$ is an indecomposable
direct
summand of $\rad P(S)$. Let $0\ne\delta\colon I(T)\to I(S)$ be a map
which is not an
isomorphism. Since $I(T)/T$ is semisimple, $\delta$ factors over
$\alpha$,
hence $g(T)=0$, or equivalently $T$ is an indecomposable direct
summand of
$Z$. Hence there is an irreducible map $f\colon
\Omega_{\what{\la}}S\to T$ in
$\uline{\mod}\what{\la}$, so there is an irreducible map $f\colon
S[-1]\to T$
in ${\bf D}^b(\la)$
\end{proof} 


\section{Components}

We consider the embedding $\mu: {\bf D}^b(\la) \to \stmodRepL$. The category
$\stmodRepL$ has almost split triangles, where for an indecomposable
$X \in \stmodRepL$, the translate $\tau_\repL X = \nu_\repL[-2] X$. In
general,
${\bf D}^b(\la)$ will not have almost split triangles. However it was shown in
\citeHII{} that for each $  P \in {\bf K}^b({_{\la}\cp})$ indecomposable there is
an
almost split triangle in ${\bf D}^b(\la)$ of the form $\nu  P[-1]
\overset{u}\to E
\overset{v}\to   P \overset{w}\to \nu  P$ where $\nu: {\bf K}^b({_{\la}\cp})
\to
{\bf K}^b({_{\la}\ci})$ is the Nakayama functor. We will show first that this triangle
is
sent under $\mu$ to the almost split triangle in $\stmodRepL$ ending
at
$\mu(  P)$ and then apply this to determine the structure of the
components
of the Auslander-Reiten quiver of ${\bf K}^b({_{\la}\cp})$ in case $\Lambda$ is a
Gorenstein
algebra.

\begin{lemma} \label{5.1}
If $  P \in {\bf K}^b({_{\la}\cp})$, then $\tau_\repL \mu(  P) \in \Img \mu$.
\end{lemma}

\begin{proof}
Let $  P = (P^i,d^i) \in {\bf K}^b({_{\la}\cp})$. Since $\mu$ commutes with the
translation functors we may assume that $P^i = 0$ for $i>0$ and it is
enough

to show that $\nu_\repL \mu(  P) \in \Img \mu$. Since $  P \in
{\bf K}^b({_{\la}\cp})$
there is $m_0$ such that $P^m = 0$ for $m<m_0$. We proceed by
induction on
$m_0$. If $m_0 = 0$, then $  P$ is a stalk complex concentrated in
degree
$0$, so $\mu(  P)$ is the stalk module $P^0$ concentrated in degree
zero.
But $\nu_\repL P^0 \cong \nu P^0[1]$ shows that $\nu_\repL P^0 \in
\Img \mu$. If
$m_0 < 0$, let $P' = (P'^i,d'^i)$ with $P'^i = P^i$ for $i<0$ and
$P'^0 = 0$,
$d'^i = d^i$ for $i<-1$ and $d'^i = 0$ for $i \ge -1$ be the truncated
complex.

We clearly have a map of complexes $P'[-1] \overset{u}\to P^0$ whose
mapping
cone is $  P$. So we obtain a triangle $P'[-1] \to P^0 \to   P
\to P'$
in ${\bf K}^b({_{\la}\cp})$. This yields a triangle $\nu_\repL \mu P'[-1] \to
\nu_\repL P^0 \to
\nu_\repL \mu   P \to \nu_\repL \mu P'$ in $\stmodRepL$.

By induction the first two terms belong to $\Img \mu$, hence so does
the
third, since $\mu$ is a triangle functor.
\end{proof}

\begin{proposition} \label{5.2}
Let $  P \in {\bf K}^b({_{\la}\cp})$ and let $\nu   P[-1] \overset{u}\to E
\overset{v}
\to   P \overset{w}\to \nu   P$ be the almost split triangle in
${\bf D}^b(\la)$
ending at $  P$. Then $\mu \nu   P[-1] \overset{\mu(u)}\to \mu E
\overset{\mu(v)}\to \mu   P \overset{\mu(w)}\to \mu \nu   P$ is
the
almost split triangle in $\stmodRepL$ ending at $\mu   P$.
\end{proposition}

\begin{proof}
Let
$$ \tau_\repL \mu   P \overset{u}\to F \overset{v}\to \mu   P
   \overset{w}\to  \tau_\repL \mu   P[1]
   \eqno{(*)}
$$
be the almost split triangle in $\stmodRepL$ ending at $\mu  
P$. By
Lemma~\ref{5.1} there is $  X \in {\bf D}^b(\la)$ such that $\mu   X =
\tau_\repL
\mu   P[1]$. So there is some $w':   P \to   X$ such that
$\mu w' =
w$. Let

$$   X[-1] \overset{u'}\to E \overset{v'}\to   P \overset{w'}\to   X
   \eqno{(**)}
$$
be a triangle in ${\bf D}^b(\la)$. By construction $\mu   X[-1] \to \mu E \to
\mu
  P \overset{w}\to \mu   X$ is isomorphic to $(*)$. Since $\mu$
is an
embedding and $(*)$ is an almost split triangle, we infer that $(**)$
is the
almost split triangle in ${\bf D}^b(\la)$ ending at $  P$.
\end{proof}

In the following let $\Lambda$ be a Gorenstein algebra. Then the
Nakayama
functor $\nu: {\bf K}^b({_{\la}\cp}) \to {\bf K}^b({_{\la}\ci})$ is an endofunctor, hence ${\bf K}^b({_{\la}\cp})$
has
almost split triangles, which are almost split triangles in ${\bf D}^b(\la)$,
and
therefore by Proposition~\ref{5.2} also almost split triangles in
$\stmodRepL$. Hence we get

\begin{corollary} \label{5.3}
Let $\C$ be a connected component of the Auslander-Reiten quiver of
${\bf K}^b({_{\la}\cp})$. Then $\C$ is a connected component of the Auslander-Reiten
quiver of $\stmodRepL$.
\end{corollary}

We will now investigate the shape of the components of the
Auslander-Reiten
quiver of ${\bf K}^b({_{\la}\cp})$ for $\Lambda$ a selfinjective algebra.

\begin{theorem} \label{5.4}
Let $\Lambda$ be a connected selfinjective algebra, which is not
semisimple.
Let $\C$ be a connected component of the Auslander-Reiten quiver of
${\bf K}^b({_{\la}\cp})$.
Then $\C$ is of the form $\ZAinf$.
\end{theorem}

\begin{proof}
Let $\C$ be a connected component of the Auslander-Reiten quiver of
${\bf K}^b({_{\la}\cp})$.
By Corollary~\ref{5.3}, $\C$ is a connected component of the
Auslander-Reiten
quiver of $\stmodRepL$. In fact $\C$ is a connected component of the
Auslander-Reiten quiver of $\modRepL$, since otherwise there would
exist an
indecomposable projective $\repL$-module $P$ such that $\rad P \in
\C$, in
particular $\rad P \in \Img \mu_{|{\bf K}^b({_{\la}\cp})}$.

But it follows from the description of $\Img \mu$ in \secIII{} that
$\Img \mu_{|{\bf K}^b({_{\la}\cp})} = \{(X_i,f_i) \mid X_i \textrm{ is
a projective } \Lambda\textrm{-module}\}$, so $\rad P \not\in \Img
\mu_{|{\bf K}^b({_{\la}\cp})}$, since $\Lambda$ is not semisimple. Consider
$l: \C \to \N$ defined by $l(X) = |X|$, the length of $X$ as a
$\repL$-module.
Then $l$ is an additive function on $\C$, since $\C$ is a component of
the
Auslander-Reiten quiver of $\modRepL$: If $  P \in {\bf K}^b({_{\la}\cp})$ and
$\mu(  P) = (X_i,f_i)$ we have that $X_i = 0$ for $|i|>m$ and some
$m$
and $X_i$ is projective for $|i| \leq m$. Since $\Lambda$ is
selfinjective
there exists $n \in \N$ such that $\nu_\Lambda^n P \cong P$ for each
projective
$\Lambda$-module $P$. If $\mu(  P) = (X_i,f_i)$, so that $l(\mu
  P) =
\sum_i |X_i|$, then $l(\Omega_\repL \mu   P) = \sum_i |\nu_\Lambda X_i|$.

So $l(\mu   P) = l(\Omega^n_\repL \mu   P)$, hence
$l(\tau^n_\repL \mu
  P) = l(\mu   P)$, showing that $l$ is a $\tau_\repL$-periodic
additive
function. Let $X \in \C$ and let $0 \ne f: P \to X$ with $P$ an
indecomposable
projective $\repL$-module. Since $\C$ does not contain any projective
$\repL$-modules we obtain for each $i$ a chain of irreducible maps
$X_i
\overset{f_i}\to X_{i-1} \to \dots \to X_1 \overset{f_1}\to X_0 = X$
such that
$f_i \dots f_1 \ne 0$ and $X_i \in \C$. By Harada-Sai (see~\citeARS{})
we know
that the length of the indecomposable modules in $\C$ is unbounded; so
$l$ is
unbounded on $\C$. In particular $\C$ contains infinitely many
$\tau_\repL$-orbits. By \citeF{} the tree class of $\C$ is $\Ainf$.
Trivially $\C$ does not contain any $\tau_\repL$-periodic vertices. So
$\C \cong \ZAinf$.
\end{proof}

If $\Lambda$ is a Gorenstein algebra, then in \secII{} we constructed
a functor $G: \stmodRepL \to {\bf D}^b(\la)$ such that $G\mu \cong
\id_{\bf D}^b(\la)$. Let
$\tilde\nu: {\bf D}^b(\la) \to {\bf D}^b(\la)$ be the equivalence induced by $\nu_\Lambda
=
D\Hom(-,{_\Lambda \Lambda})$, then we obtain a commutative diagram
$$
\begin{CD}
\stmodRepL  @>G>>   {\bf D}^b(\la)            \\
@VV{\nu_\repL}V     @VV{\tilde\nu[1]}V  \\
\stmodRepL  @>>>    {\bf D}^b(\la)
\end{CD}
$$

\begin{proposition} \label{5.5}
Let $\Lambda$ be a Gorenstein algebra and $X,Y \in {\bf D}^b(\la)$
indecomposable.
If $f: X \to Y$ is irreducible and $Y \not\in {\bf K}^b({_{\la}\cp})$, then
$X \cong \tilde\nu Y[-1]$.
\end{proposition}

\begin{proof}
We consider the almost split triangle in $\stmodRepL$
$$ \tau_\repL \mu Y \overset{\alpha}\to E \overset{\beta}\to \mu Y
   \overset{\gamma}\to \tau_\repL \mu Y[1]
   \eqno{(*)}
$$
Since $f: X \to Y$ is irreducible and by Theorem 4.3, also
$\mu(f)$ is
irreducible, we see that $E \cong \mu(X) \oplus C$ and $\beta =
(\mu(f),g)^t$
for some $g: C \to \mu Y$. Since $\tau_\repL = \nu_\repL
\Omega^2_\repL$ and using
the diagram above we see that $G$ applied to $(*)$ yields a triangle in ${\bf D}^b(\la)$
$$ \tilde\nu Y[-1] \overset{G(\alpha)}\to X \oplus GC
   \overset{G(\beta)}\to Y \overset{G(\gamma)}\to \tilde\nu Y
   \eqno{(**)}
$$
We claim that $G(\gamma) = 0$. Otherwise, let $h: Z \to Y$ be a map
which is
not a split epi, then $\mu(h)$ is not a split epi. But then
$\mu(h)\gamma = 0$,
 hence $0 = G(\mu(h)\gamma) = h G(\gamma)$. Since $Y$ and $\tilde\nu
 Y[-1]$
are indecomposable, $(**)$ would be an almost split triangle in
${\bf D}^b(\la)$.
Since $Y \not\in {\bf K}^b({_{\la}\cp})$ this contradicts \citeHII{}. So $G(\gamma) =
0$,
hence $G(\beta)$ is a split epi. Since $X$ is not isomorphic to $Y$,
we get
that $X \cong \tilde\nu Y[-1]$.
\end{proof}

We will now show that for selfinjective algebras $\Lambda$ irreducible
maps
in ${\bf D}^b(\la)$ outside ${\bf K}^b({_{\la}\cp})$ are rare. For this we will need the
following
easy fact, but first we will define the relevant class of algebras
$\Lambda_n$ for $n \ge 1$. Let $\Lambda_1 = k[x]/(x^2)$ and
$\Lambda_n$ for
$n \ge 2$ defined by the following quiver
$$
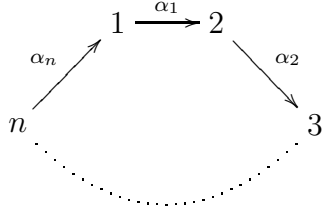
\xymatrix{
& 1 \ar[r]^{\alpha_1}   & 2 \ar[rd]^{\alpha_2}  &   \\
n \ar[ur]^{\alpha_n} \ar@{..}@/_2.5pc/[rrr] &&& 3    \\
\\
}
$$
with relations $\alpha_i \alpha_{i+1} = 0 \textrm{ for } 1 \leq i \leq
n$ where
$\alpha_{n+1} = \alpha_1$.

We collect the relevant information in the following well-known lemma
(see~\citeARS{})

\begin{lemma} \label{5.6}
Let $\Lambda$ be a basic selfinjective algebra which is not semisimple.
\begin{enumerate}
\item $\rad^2 \Lambda = 0$ if and only if $\Lambda \cong \Lambda_n$
  for some
$n$.
\item If $S(i)$ is a simple $\Lambda_n$-module, then $\nu S(i) =
  S(i-1)$
where $S(0) = S(n)$.
\item $\Lambda_n$ is symmetric if and only if $n = 1$.
\end{enumerate}
\end{lemma}

\begin{theorem} \label{5.7}
Let $\Lambda$ be a basic selfinjective algebra which is not
semisimple.
Let $Y \in {\bf D}^b(\la) \setminus {\bf K}^b({_{\la}\cp})$ be indecomposable. There exists an
irreducible map $f: X \to Y$ in ${\bf D}^b(\la)$ if and only if $\Lambda \cong
\Lambda_n$,
 $Y \cong S(i-1)[j]$, $X \cong S(i)[j-1]$ for some $1 \leq i \leq n$
 and $j \in
\Z$.
\end{theorem}

\begin{proof}
If $\Lambda = \Lambda_n$ for some $n$, we have seen in
Proposition 4.4
that for each arrow $\alpha_i$ we have an irreducible map $\nu
S(i)[-1] \to
S(i+1)$ in ${\bf D}^b(\la)$. Since $\pd_{\Lambda_n} S(i+1) = \infty$ we have $S(i+1) \not\in {\bf K}^b({_{\la}\cp})$.

Conversely, let $f: X \to Y$ be irreducible in ${\bf D}^b(\la)$ and $Y \not\in
{\bf K}^b({_{\la}\cp})$.
We choose $Y \in {\bf K}^{-,b}({_{\la}\cp})$ and may assume that $Y = (P^i,d^i)$
satisfies $Y^i=0$
 for $i>0$ and $H^0(Y) \ne 0$. By Proposition~\ref{5.5} we know that
 $X \cong
\tilde\nu Y[-1]$. Since $\Lambda$ is selfinjective, we have that
$\nu_\Lambda$
is exact, hence $\tilde\nu Y = (\nu_\Lambda P^i, \nu_\Lambda
d^i)$. Consider the
triangle $\nu_\Lambda P^0[-1] \overset{\alpha}\to Y[-1]
\overset{\beta}\to
Y[-1]_{\leq 0} \overset{\gamma}\to \nu_\Lambda P^0$. Since
$\Hom(\nu_\Lambda P^0[-1],Y)
 = 0$, there is some $\bar f$ such that $f = \beta \bar f$. $H^0(Y)
 \ne 0$
implies that $\alpha \ne 0$. Thus $\beta$ is not split mono, so $\bar
f$ is
split epi, since $f$ is irreducible. Hence $H^i(\bar f): H^i(\tilde\nu
Y[-1]_{\leq 0} \to H^i(Y)$ is split epi for all $i$. Since $\nu$ is
exact,
we have that $H^i(\tilde\nu Y) \cong \nu_\Lambda H^i(Y)$ for all
$i$. Also we
have that $H^i(\tilde\nu Y[-1]_{\leq 0}) = H^i(\tilde\nu Y[-1])$ for
all $i
\leq -1$ and $H^i(\tilde\nu Y[-1]) = H^{i-1}(\tilde\nu Y)$. Since $Y
\in {\bf K}^{-,b}({_{\la}\cp})$ there is $n_0 \leq 0$ such that $H^n(Y) = 0$ for all $n \leq
n_0$.
Choose $n_0$ maximal with this property, so $H^{n_0}(Y) = 0$ and
$H^{n_0+1}(Y)
\ne 0$. We claim that $Y \in \modL$ or equivalently $n_0 =
-1$. Otherwise
$n_0 \leq -2$. But then $0 = H^{n_0}(Y) = H^{n_0}(\tilde\nu Y) =
H^{n_0+1}(\tilde\nu Y[-1]) = H^{n_0+1}(\tilde\nu Y[-1]^{\leq 0})
\twoheadrightarrow H^{n_0+1}(Y)$, hence $H^{n_0+1}(Y) = 0$ in contrast
to the
choice of $n_0$. Hence $n_0 = -1$, so $Y$ is an indecomposable
$\Lambda$-module.
 But then $Y^{\leq-1} \cong \Omega_\Lambda Y[1]$, so $Y^{\leq -1}$ is
 indecomposable,
since $\Lambda$ is selfinjective. But then $\tilde\nu Y[-1]^{\leq 0}$
is
indecomposable, so $\bar f$ is an isomorphism, so $\nu_\Lambda
\Omega_\Lambda Y
\cong Y$. If $Y$ is not simple, there is a proper epi $Y
\overset{\pi}\to S$
for some simple $S$. So there is $h: Y \to I(S)$, with $\Ker h \ne 0$
and
$\Coker h \ne 0$. Since $\Hom(I(S),Y[1]) = \Ext^1_\Lambda(I(S),Y) = 0$
we obtain
a triangle
$$ C_h[-1] \overset{g}\to Y \overset{h}\to I(S) \to C_h $$
with $C_h[-1]$ indecomposable and $H^0(C_h[-1]) = \Ker h$,
$H^1(C_h[-1]) =
\Coker h$. Since $\Hom(\nu_\Lambda Y[-1],I(S)) = 0$, there is $f':
\nu_\Lambda
Y[-1] \to C_h[-1]$ such that $f'g = f$. Since $h \ne 0$, $g$ is not
split epi,
hence $f'$ is split mono, since $f$ is irreducible. Since $C_h[-1]$ is
indecomposable, we have that $f'$ is an isomorphism, in contrast to
$H^0(C_h[-1]) \ne 0 \ne H^1(C_h[-1])$, so $X$ is a simple
$\Lambda$-module.
Since $Y \cong \nu_\Lambda \Omega_\Lambda Y$, we see that
$\Omega_\Lambda Y$ is a
simple $\Lambda$-module. But then $\rad^2 \Lambda = 0$, since
$\Lambda$ is
selfinjective and the assertion follows from Lemma~\ref{5.6}.
\end{proof}

\section{Behavior of irreducible maps}

In this section we show that beyond the Gorenstein algebras the
behavior of
irreducible maps in ${\bf D}^b(\la)$ is not so regular. In particular, we show
that
some natural conjectures have a negative answer.

For a non-zero map $f:   P \to   Q$ between indecomposable
objects in
${\bf D}^b(\la)$ but not in ${\bf K}^b({_{\la}\cp})$ we investigate the connection between $f:
  P
\to   Q$ and $f_{\geq -n}:   P_{\geq -n} \to   Q_{\geq -n}$
being irreducible
for some $n$. We also give some sufficient condition for an
irreducible map
in $\modL$ not be irreducible in ${\bf D}^b(\la)$.

We start with a general result on mapping cones of irreducible maps,
where
the analogous result in abelian categories is well known.

\begin{proposition} \label{6.1}
Let $X$ and $Y$ be indecomposable in ${\bf D}^b(\la)$, for a finite dimensional
algebra $\Lambda$, and assume that we have an irreducible map $f: X
\to Y$.
Then the mapping cone $C_f$ is indecomposable.
\end{proposition}

\begin{proof}
This can be proved in a similar way as the abelian analog. Here we
give a
slightly shorter proof using Theorem 4.3.
Let $\mu: {\bf D}^b(\la) \to \stmodRepL$ be as usual the
natural embedding. Then we know from Theorem 4.3 that $\mu(f):
\mu(X)
\to \mu(Y)$ is irreducible in $\stmodRepL$. This is induced by an
irreducible
map $f': \mu(X) \to \mu(Y)$ in $\modRepL$. If $f'$ is mono, we have an
exact
sequence $0 \to \mu(X) \to \mu(Y) \to \mu(Y)/\mu(X) \to 0$, and if
$f'$ is
epi, we have an exact sequence $0 \to \Ker f' \to \mu(X) \to \mu(Y)
\to 0$.
We know that in the first case $\mu(Y)/\mu(X)$ is indecomposable and
in the
second case $\Ker f'$ is indecomposable, see \cite{ARS}. So in any
case we have a triangle
$\mu(X) \overset{\mu(f)}\to \mu(Y) \to Z$ in $\stmodRepL$, where $Z$
is
indecomposable. Since $\mu(C_f) \cong Z$, it follows that $C_f$ is
indecomposable.
\end{proof}

Let $0 \to A \to B \to C \to 0$ be an almost split sequence in
$\modL$. Then
it is known that if $\idim_\Lambda A \leq 1$ and $\pd_\Lambda C \leq 1$,
then the
sequence gives rise to an  almost split triangle in ${\bf D}^b(\la)$ \citeHI.
Consequently the
corresponding
irreducible maps $f_i: A \to B_i$ and $g_i: B_i \to C$ stay
irreducible, where
$B = \bigoplus_{i=1}^t B_i$ with $B_i$ indecomposable. But the normal
behavior
is that irreducible maps in $\modL$ do not stay irreducible in
${\bf D}^b(\la)$. We
illustrate this with the following result.

\begin{proposition} \label{6.2}
Let $\Lambda$ be a finite dimensional algebra, $X$ and $Y$
indecomposable
$\Lambda$-modules with $\pd_\Lambda X < \infty$ and $\pd_\Lambda Y \ge
\pd_\Lambda X + 2$. Then there is no irreducible map $f: X \to Y$ in ${\bf D}^b(\la)$.
\end{proposition}

\begin{proof}
Assume that we have an irreducible map $f: X \to Y$ in $\modL$, with
$X$ and
$Y$ indecomposable, $\pd X = i < \infty$ and $\pd Y \geq i + 2$. Let
$  P:
\dots \to P^{-(i+2)} \to \dots \to P^{-1} \to P^0 \to 0$ be a minimal
projective
resolution of $Y$ and $  Q: 0 \to Q^{-i} \to \dots \to Q^{-1} \to
Q^0 \to 0$
a minimal projective resolution of $X$. Let $  C$ denote the
complex $0
\to P^{-(i+1)} \to \dots \to P^{-1} \to P^0 \to 0$. Then $f:   Q
\to   P$
factors as $  Q \overset{  g}\to   C \overset{  h}\to  
P$,
since we have the commutative diagram
\smallskip
{\minCDarrowwidth1.5pc
$$
\begin{CD}
    @.         @.      0       @>>> Q^{-i}  @>>> \dots
@>>> Q^{-1} @>>> Q^0  @>>>  0                       \\
  @.            @.       @VVV       @V{f^{-i}}VV    @.
   @V{f^{-1}}VV   @V{f^0}VV                     \\
    @.      0      @>>> P^{-(i+1)} @>>> P^{-i}  @>>> \dots
@>>> P^{-1} @>>> P^0  @>>>  0                   \\
  @.            @.        @|           @|       @.
    @|        @|                        \\
\dots @>>>  P^{-(i+2)} @>>> P^{-(i+1)} @>>> P^{-i}  @>>> \dots
@>>> P^{-1} @>>> P^0  @>>>  0                   \\
\end{CD}
$$
}
\smallskip

We want to show that $  g$ is not a split monomorphism and $  h$
is
not a split epimorphism. If $  g:   Q \to   C$ was a split
monomorphism, the induced map $H^0(  Q) = X \to H^0(  C) = Y$
would
be a split monomorphism. Since $X$ and $Y$ are indecomposable
nonisomorphic
modules, this is impossible.

The diagram
\smallskip
{\minCDarrowwidth1.5pc
$$
\begin{CD}
\Omega^{-(i+2)} Y                   \\
    @VVV                        \\
  P^{-(i+1)}    @>>> \dots @>>> P^{-1}   @>>>   P^0 \\
    @.          @.   @.     @VVV    \\
        @.     @.        @.  Y  \\
\end{CD}
$$
}
\smallskip

\noindent
gives rise to the triangle $\Omega^{-(i+2)} Y[i+1] \to   C \to Y[0]
\overset{\alpha}\to \Omega^{-(i+2)} Y[i+2]$. Here $\alpha$ is given by
the
sequence $0 \to \Omega^{-(i+2)} Y \to P^{-(i+1)} \to \dots \to P^{-1}
\to P^0 \to Y
\to 0$, which does not represent the zero element since $\pd_\Lambda Y
\ge i + 2$. Hence, $  h:   C \to Y[0]$ or equivalently $  h:
  C
\to   P$ is not a split epimorphism.
It follows that $  f:   Q \to P$ is not irreducible.
\end{proof}

Note that if $\Lambda$ is hereditary, then each irreducible map in
$\modL$
stays irreducible in ${\bf D}^b(\la)$. In this case $\pd X$ is $0$ or $1$, hence
we
can never have $\pd Y \geq \pd X + 2$.

The next natural question is to which extent we have irreducible maps
$X \to Y[1]$, where $X$ and $Y$ are indecomposable in $\modL$,
corresponding
to elements of $\Ext^1_\Lambda(X,Y)$. Here we have seen some
sufficient
conditions in \secIII{}. Normally, we do not have such irreducible maps.

\begin{proposition} \label{6.3}
Let $f: X \to Y[1]$ be an irreducible map, where $X$ and $Y$ are
indecomposable $\Lambda$-modules. Then $Y$ must be a summand of $\Omega X$.
\end{proposition}

\begin{proof}
We have the factorization $X \overset{h}\to \Omega X[1]
\overset{g[1]}\to Y[1]$ of $f: X \to Y[1]$, as is seen by considering
the
diagram
$$
\xymatrix{
0 \ar[r]            &
\Omega X \ar[r] \ar@{.>}[d]^g   &
P_X \ar[r] \ar@{.>}[d]      &
X \ar[r] \ar@{=}[d]     &
0
\\
0 \ar[r]            &
Y \ar[r]            &
E \ar[r]            &
X \ar[r]            &
0
}
$$

Then $h: X \to \Omega X[1]$ is not a split monomorphism since $H^0(h)$
is
not a split monomorphism. Since $f: X \to Y[1]$ is irreducible, it
follows
that $g: \Omega X \to Y$ is a split epimorphism, so that $Y$ is a
summand
of $\Omega X$.
\end{proof}

We now give another situation where there are no irreducible maps,
containing the case $X\to Y[2]$, corresponding to elements of
$\Ext^2_\Lambda(X,Y)$, as a special case.

\begin{proposition} \label{6.4}
Let $ P$ and $ Q$ be indecomposable objects in ${\bf D}^b(\la)$ for a
finite dimensional algebra $\Lambda$, represented by complexes of
projective $\Lambda$-modules with no split exact summands, with
$P^0 \ne 0$, $P^i = 0$ for $i>0$ and $Q^i = 0$ for $i \ge -1$. Then
there is
no irreducible map $f:  P \to  Q$ in ${\bf D}^b(\la)$.
\end{proposition}

\begin{proof}
Let $f:  P \to  Q$ be a map in ${\bf D}^b(\la)$. Consider the
factorization
of $f$ given by
\smallskip
{\minCDarrowwidth1.5pc
$$
\begin{CD}
\dots @>>>    P^{-2} @>>>     P^{-1} @>a>>    P^0 @>>>   0   @.
\qquad\qquad @.  P                          \\
  @.        @VVV        @VVV          @VVV  @.
@.      @VVhV                           \\
\dots @>>>    P^{-2} @>>>     P^{-1} @>>>     0   @>>>   0   @.
\qquad\qquad @.  P_{\leq -1}                    \\
  @.        @VVV        @VVV          @.    @.
@.      @VVgV                           \\
\dots @>>>    Q^{-2} @>>>   0    @.       @.         @.
\qquad\qquad @. Q                          \\
\end{CD}
$$
}
\smallskip

\noindent
We have $H^0(  P) = P^0/\Img a$, which is not zero since $  P$
has no
split exact direct summands. Since $H^0(  P_{\leq -1}) = 0$, $h:
  P \to
  P_{\leq -1}$ cannot be a split monomorphism.

Assume now that $g$ is a split epimorphism, and consider the triangle
$  P_{\leq -1} \overset{g}\to   Q \overset{u}\to C_g \to$. Then
$u:   Q \to C_g$ must be homotopic to $0$, that is, we have the diagram
$$
\xymatrix{
\dots \ar[r]                            &
Q^{-3} \ar[dl] \ar[d]^{(0,1)} \ar[r]^{b_{-2}}           &
Q^{-2} \ar[dl]^{s_{-2}} \ar[d]^{(0,1)} \ar[r]^{b_{-1}}      &
\quad 0 \quad \ar[dl]^{s_{-1}}
\\
\dots \ar[r]                            &
P^{-2} \oplus Q^{-3} \ar[r]_{c_{-2}}                &
P^{-1} \oplus Q^{-2} \ar[r]_{c_{-1}}                &
\quad 0 \quad
}
$$
where  $b_{-i} s_{-i}  +s_{-(i+1)} c_{-(i+1)}  = (0,1)$ for all $i \geq
1$. Using
the same maps $s_i$ we see that in the triangle $  P \overset{f}\to
  Q \overset{v}\to C_f \to$, the map $v$ must be $0$, so that $f$
would
also be a split epimorphism. Since $  P$ and $  Q$ are
indecomposable,
$f$ would be an isomorphism, which is impossible because $H^0(  P)
\ne 0$
and $H^0(  Q) = 0$. We conclude that $g$ is not a split
epimorphism.
Since we already have that $h$ is not a split monomorphism, it follows
that $f:   P \to   Q$ is not irreducible.
\end{proof}

The following sufficient condition for the mapping cone to be
indecomposable
will be useful. We formulate the result in the more general setting of
Krull-Schmidt triangulated categories.

\begin{lemma} \label{6.5}
Let $f: P \to Q$ be a map between indecomposable objects in a
Krull-Schmidt
triangulated category $\C$ with shift $[1]$ and assume that $f$ is not
zero and not invertible. Complete to a triangle
$P \overset{f}\to Q \overset{g}\to C \to P[1]$. If $\Hom(Q,P[1]) = 0$,
then $C$ is indecomposable.
\end{lemma}

\begin{proof}
Assume to the contrary that $\C$ is not indecomposable, and write
$C = \bigoplus_{i=1}^r Z_i$, where $r > 1$ and each $Z_i$ is
indecomposable.
Let $g = (g_1, \dots, g_r)$ and $h = (h_1, \dots, h_r)^t$. Then we
know from
\citeR{} that $g_i \ne 0$ and $h_i \ne 0$ for each $i = 1, \dots, r$.

Consider the map
$$
\varphi = \left(
\begin{array}{cccc}
     1 & 0     & \dots & 0      \\
     0 &       &       & \vdots     \\
\vdots &       &       & \vdots     \\
     0 & \dots & \dots & 0
\end{array}
\right)
: C \to C,
$$
where $1 = \id_{Z_1}$. We then have the diagram
$$
\begin{CD}
P   @>f>>   Q   @>g>>   C   @>h>>   P[1]        \\
@.      @.      @VV{\varphi}V           \\
P   @>f>>   Q   @>g>>   C   @>h>>   P[1]        \\
\end{CD}
$$
Since by the assumption $\Hom(Q,P[1]) = 0$, it follows that $g\varphi
h = 0$.
Hence there is a map $\varphi_Q: Q \to Q$ such that $\varphi_Q g = g
\varphi$.
We have $g\varphi = (g_1, 0, \dots, 0)$ and $\varphi_Q g = (\varphi_Q
g_1,
\dots, \varphi_Q g_r)$, so that $\varphi_Q g_1 = g_1$ and $\varphi_Q
g_i = 0$
for $2 \leq i \leq r$. Since $\varphi^2 = \varphi$, we have
$\varphi_Q^n g = g\varphi^n = g\varphi$, so that $\varphi_Q^n g_1 =
g_1$ and
$\varphi_Q^n g_i = 0$ for $2 \leq i \leq r$. Since $Q$ is
indecomposable,
any map $t: Q \to Q$ is nilpotent or an isomorphism, so that we have a
contradiction. It follows that $C$ is indecomposable.
\end{proof}

We now consider the following question. If we have an irreducible map
$ f: P \to  Q$ between unbounded complexes of projective
modules, not objects in ${\bf K}^b({_{\la}\cp})$, is then $  f_{\ge -n}:
P_{\ge -n}
\to  Q_{\ge -n}$ irreducible for all $n$ where $ f_{\ge -n}$
is a
nonzero map between indecomposable objects?

For selfinjective algebras $\Lambda$, the existence of an irreducible
map $ f:  P \to  Q$ not in ${\bf K}^b({_{\la}\cp})$ implies that
$\Lambda$ is
selfinjective with $\rad^2 \Lambda = 0$ and that we have $  f: S
\to T[1]$,
where $S$ and $T$ are simple $\Lambda$-modules. In this case $
f_{\ge -n}:
  P_{\ge -n} \to   Q_{\ge -n}$ is irreducible for $n \ge 2$.

We now give an example which gives a negative answer to the above
question.
Let $\Lambda$ be the path algebra of the quiver
\smallskip
$$
\xymatrix{
1 \ar@(ul,dl)[]_{\gamma} &
2 \ar[l]_{\alpha}
}
$$
with relations $\gamma\alpha = 0$, $\gamma^2 = 0$. Denote by $S$ the
simple
module at vertex $1$ and by $T$ the simple module at vertex $2$. Then
the
indecomposable projective $\Lambda$-modules have Loewy series ${^S_S}$
and
${^T_S}$, and we have $T$ and ${^{ST}_{\,S}}$ for the indecomposable
injectives.
We know from Proposition 4.4 that the map $f: T \to S[1]$ is
irreducible, and
we can write this as $  f:   P \to   Q$ given by
$$
\begin{CD}
\dots @>>>  {^S_S} @>>> {^S_S} @>>> {^T_S}      \\
@.      @|      @|      @VVV        \\
\dots @>>>  {^S_S} @>>> {^S_S} @>>> 0       \\
\end{CD}
$$

We then have the following.

\begin{proposition} \label{6.6}
Let $\Lambda$ be as above. In the above notation we have that
$  f:   P \to   Q$ is irreducible, while $  f_{\ge -1}:
  P_{\ge -1} \to   Q_{\ge -1}$ is a map between indecomposable
objects
which is not irreducible.
\end{proposition}

\begin{proof}
We have already seen that $  f:   P \to   Q$ is
irreducible. We now
want to show that $  f_{\ge -1}:   P_{\ge -1} \to   Q_{\ge
  -1}$ is not
irreducible. We have $  Q_{\ge -1} = ({^S_S})[1]$, and hence
$\nu   Q_{\ge -1} = {^{ST}_{\,S}}$, so that we have an almost split
triangle ${^{ST}_{\,S}} \to ({^S_S} \overset{\alpha}\to {^{ST}_{\,S}})
\to {^S_S}[1] \to$, where $\Img \alpha = S$. We claim that $X =
({^S_S}
\overset{\alpha}\to {^{ST}_{\,S}})$ is indecomposable. For this, it is
sufficient to show that $\Hom({^{ST}_{\,S}}, {^S_S}[1]) = 0$ by
Lemma~\ref{6.5},
that is that $\Ext^1_\Lambda({^{ST}_{\,S}}, {^S_S}) = 0$. This follows
by
considering the injective resolution $0 \to {^S_S} \overset{h}\to
{^{ST}_{\,S}}
\overset{g}\to T \to 0$, which gives rise to the exact sequence
$\Hom({^{ST}_{\,S}}, {^{ST}_{\,S}}) \overset{\varphi}\to
\Hom({^{ST}_{\,S}},T) \to
\Ext^1({^{ST}_{\,S}}, {^S_S}) \to 0$, and using that $\varphi$ is
clearly an
epimorphism. Hence we conclude that $X$ is
indecomposable. Alternatively
we could prove that $X$ is indecomposable by considering the
homology
of $X$ and how it could decompose.

Since $H^0(  P_{\ge -1}) = T$ while $H^0(X) = S \oplus T$, $
P_{\ge -1}$
cannot be isomorphic to $X$. Hence $  f_{\ge -1}:   P_{\ge -1}
\to
  Q_{\ge -1}$ is not irreducible.
\end{proof}

We now give an example of a nonzero map $  f:   P \to   Q$
between
indecomposable objects which is not irreducible, but such that
$  f_{\ge -n}:   P_{\ge -n} \to   Q_{\ge -n}$ is an
irreducible map
between indecomposable objects for some $n$.

Let $\Lambda = k[x]/(x^3)$, and consider the complexes of projective modules:
\begin{align*}
  P \qquad\qquad
&\dots \to
\Lambda \overset{\cdot x}\to
\Lambda \overset{\cdot x^2}\to
\Lambda \overset{(x,0)}\to
\Lambda \oplus \Lambda \overset{({^{x^2}_{\,x}})}\to
\Lambda
\\
  Q \qquad\qquad
&\dots \to
\Lambda \overset{\cdot x}\to
\Lambda \overset{\cdot x^2}\to
\Lambda \; \overset{\cdot x}\to \quad\;\,
\Lambda \quad\;\, \to
0
\end{align*}
where the right hand terms are in degree $0$, as objects in ${\bf D}^b(\la)$.
Consider the map $  f:   P \to   Q$ in ${\bf D}^b(\la)$ induced by the
commutative diagram
$$
\begin{CD}
\dots @>>>
\Lambda @>{\cdot x}>>
\Lambda @>{\cdot x^2}>>
\Lambda @>{(x,0)}>>
\Lambda \oplus \Lambda @>{({^{x^2}_{\,x}})}>>
\Lambda
\\
@. @| @| @| @VV{({^1_0})}V @VVV
\\
\dots @>>>
\Lambda @>{\cdot x}>>
\Lambda @>{\cdot x^2}>>
\Lambda @>{\cdot x}>>
\Lambda @>>>
0
\end{CD}
$$

We have the following

\begin{proposition} \label{6.7}
With the above notation and assumptions we have the following
\begin{enumerate}
\item The induced map $  f_{\ge -2}:   P_{\ge -2} \to
  Q_{\ge -2}$ is
an irreducible map between indecomposable objects in ${\bf D}^b(\la) $.
\item The map $  f:   P \to   Q$ is a map between
  indecomposable
objects which is not irreducible.
\end{enumerate}
\end{proposition}

\begin{proof}
(1) The map $  f_{\ge -2}$ is given by the following diagram
$$
\begin{CD}
0 @>>>
\Lambda @>{(x,0)}>>
\Lambda \oplus \Lambda @>{({^{x^2}_{\,x}})}>>
\Lambda @>>>
0
\\
@. @| @VV{({^1_0})}V @VVV
\\
0 @>>>
\Lambda @>{\cdot x}>>
\Lambda @>>>
0 @>>>
0
\end{CD}
$$

Since $\Lambda = k[x]/(x^3)$ is symmetric, we have $\tau(  Q_{\ge
  -2}) =
  Q_{\ge -2}[-1]$, and hence an almost split triangle $  Q_{\ge
  -2}[-1]
\to E \to   Q_{\ge -2} \overset{\alpha}\to   Q_{\ge -2}$. The
map
$\alpha:   Q_{\ge -2} \to   Q_{\ge -2}$ inducing the almost
split
triangle is easily seen to be given by the diagram
$$
\begin{CD}
0   @>>>    \Lambda @>{\cdot x}>>   \Lambda @>>>    0   \\
@.      @VV{0}V         @V{\cdot x^2}VV     \\
0   @>>>    \Lambda @>{\cdot x}>>   \Lambda @>>>    0   \\
\end{CD}
$$
For it is clear that the induced map is nonzero and is in the socle of
$\End(  Q_{\ge -2})$. Taking the mapping cone of $\alpha$ we obtain
$  P_{\ge -2}[1]$, so that $E \cong   P_{\ge -2}$. This shows
that
$  f_{\ge -2}:   P_{\ge -2} \to   Q_{\ge -2}$ is irreducible.

We next show that $  P_{\ge -2}$ is indecomposable. We give
a
proof which at the same time illustrates the previous theory, rather
than
giving a direct computational proof.
We know from Theorem 5.4 that the components of the AR-quiver of
${\bf K}^b({_{\la}\cp})$ are of the
form
$\ZAinf$, and that the image of a component for ${\bf K}^b({_{\la}\cp})$ is a
component
of the AR-quiver for $\stmodRepL$. All $\repL$-modules in such a
component $\C$ are given by projective modules, the same ones as for
${\bf K}^b({_{\la}\cp})$. Then $\C$ is also a component for $\modRepL$. This follows
since any indecomposable projective object in $\modRepL$ has an
irreducible
map to this object modulo its socle, and this object is not given by
only
projective modules.

If $  P_{\ge -2}$ was not indecomposable, then
$\Lambda \overset{\cdot x}\to \Lambda$ would not be at the border of
the $\ZAinf$-component. Hence we would have an irreducible epimorphism
starting at $\Lambda \overset{\cdot x}\to \Lambda$, which would then
have
to end at $\Lambda$, since the terms must be projective. But on the
other
hand we have an almost split triangle $\Lambda[-1] \to
(\Lambda \overset{\cdot x^2}\to \Lambda) \to \Lambda$ in ${\bf K}^b({_{\la}\cp})$,
which
gives a contradiction.

(2) We first show that $  P$ and $  Q$ are indecomposable. This
is
obvious for $  Q$. Assume $  P =  {P'} \oplus  {P''}$ is a
nontrivial decomposition. Then we have $  P_{\ge-2} =
 {P'}_{\ge-2}
\oplus  {P''}_{\ge-2}$. Since $H^0(  P) = S$, $H^{-1}(  P) =
{^S_S}$
and $H^{-i}(  P) = 0$ for $i \ne 0,1$, we must have, say $ {P'}
\cong S$
and $ {P''} \cong ({^S_S})[1]$. But then $ {P'}_{\ge-2}$ and
$ {P''}_{\ge-2}$ are both nonzero, contradicting that $
P_{\ge-2}$ is
indecomposable.

That $  f:   P \to   Q$ is not irreducible follows since
$\Lambda$
is selfinjective and $\rad^2 \Lambda \ne 0$ and $  P$ and $  Q$
are
not in ${\bf K}^b({_{\la}\cp}) $. We could alternatively give a direct argument by
considering the following factorization of the map $  f:   P \to   Q$:
$$
\begin{CD}
\dots @>>>
\Lambda @>{\cdot x}>>
\Lambda @>{\cdot x^2}>>
\Lambda @>{(x,0)}>>
\Lambda \oplus \Lambda @>{({^{x^2}_{\,x}})}>>
\Lambda @.
\qquad\qquad @.
  P
\\
@. @| @| @VV{(1,0)}V @| @| @. @VVgV
\\
\dots @>>>
\Lambda @>{\cdot x}>>
\Lambda @>{(x^2,0)}>>
\Lambda \oplus \Lambda @>{({^x_0} {^{\,0}_{x^2}})}>>
\Lambda \oplus \Lambda @>{({^{x^2}_{\,x}})}>>
\Lambda @.
\qquad\qquad @.
  U
\\
@. @| @| @VV{({^1_0})}V @VV{({^1_0})}V @VVV @. @VVhV
\\
\dots @>>>
\Lambda @>{\cdot x}>>
\Lambda @>{\cdot x^2}>>
\Lambda @>{\cdot x}>>
\Lambda @>>>
0 @.
\qquad\qquad @.
  Q
\end{CD}
$$
and showing that $g$ is not a split monomorphism and $h$ is not a
split
epimorphism. The first claim follows directly by considering the
homology of
$  P$ and $ U$, and the second claim is also not hard to show.
\end{proof} 

\end{document}